\documentclass[11pt,a4paper]{article}

\usepackage{a4wide,amssymb,graphics,graphicx,textcomp}
\usepackage{enumerate,textcomp,multirow,amsmath}
\usepackage{algorithm}
\usepackage{algorithmic}
\usepackage{url}

\newtheorem{theorem}{Theorem}[section]
\newtheorem{lemma}[theorem]{Lemma}
\newtheorem{remark}[theorem]{Remark}
\newtheorem{definition}[theorem]{Definition}
\newtheorem{corollary}[theorem]{Corollary}

\newtheorem{example}[theorem]{Example}

\numberwithin{equation}{section}
\numberwithin{table}{section}
\numberwithin{figure}{section}
\newcommand{\proof}{\par\noindent{\bf Proof}. \ignorespaces}
\newcommand{\eproof}{\space
    {\ \vbox{\hrule\hbox{\vrule height1.3ex\hskip0.8ex\vrule}\hrule}}\par}

\newcommand {\mat}  [1] {\left[\begin{array}{#1}}
\newcommand {\rix}      {\end{array}\right]}

\catcode`@=11     
\font\tenex=cmex10 
\newdimen\p@renwd
\setbox0=\hbox{\tenex B} \p@renwd=\wd0 
\def\bmat#1{\begingroup \m@th
  \setbox\z@\vbox{\def\cr{\crcr\noalign{\kern2\p@\global\let\cr\endline}}%
    \ialign{$##$\hfil\kern2\p@\kern\p@renwd&\thinspace\hfil$##$\hfil
      &&\quad\hfil$##$\hfil\crcr
      \omit\strut\hfil\crcr\noalign{\kern-\baselineskip}%
      #1\crcr\omit\strut\cr}}%
  \setbox\tw@\vbox{\unvcopy\z@\global\setbox\@ne\lastbox}%
  \setbox\tw@\hbox{\unhbox\@ne\unskip\global\setbox\@ne\lastbox}%
  \setbox\tw@\hbox{$\kern\wd\@ne\kern-\p@renwd\left[\kern-\wd\@ne
    \global\setbox\@ne\vbox{\box\@ne\kern2\p@}%
    \vcenter{\kern-\ht\@ne\unvbox\z@\kern-\baselineskip}\,\right]$}%
  \null\;\vbox{\kern\ht\@ne\box\tw@}\endgroup}
\catcode`@=12    

\newcommand{\backmatter}{\appendix
\def\chaptermark##1{\markboth{%
\ifnum  \c@secnumdepth > \m@ne  \@chapapp\ \thechapter:  \fi  ##1}{%
\ifnum  \c@secnumdepth > \m@ne  \@chapapp\ \thechapter:  \fi  ##1}}%
\def\sectionmark##1{\relax}}

\makeatletter
\newcommand*{\rom}[1]{\expandafter\@slowromancap\romannumeral #1@}
\makeatother

\def\real{\mathop{\mathrm{Re}}}
\def\imag{\mathop{\mathrm{Im}}}

\newif\ifMDlatex
\catcode`@=11     
\def\MD@us#1{\csname#1style\endcsname}
\def\MD@uf#1{\csname#1font\endcsname}
\def\MD@t{text}
\def\MD@s{script}
\def\MD@ss{scriptscript}
\newdimen\MD@unit
\MD@unit\p@
\def\MD@changestyle#1{
  \relax\MD@unit0.1\fontdimen6\MD@uf{#1}0
  \everymath\expandafter{\the\everymath\MD@us{#1}}
}
\def\MD@dot{$\m@th\ldotp$}
\def\MD@palette#1{\mathchoice{#1\MD@t}{#1\MD@t}{#1\MD@s}{#1\MD@ss}}
\def\MD@ddots#1{{\MD@changestyle{#1}%
  \mkern1mu\raise7\MD@unit\vbox{\kern7\MD@unit\hbox{\MD@dot}}%
  \mkern2mu\raise4\MD@unit\hbox{\MD@dot}%
  \mkern2mu\raise \MD@unit\hbox{\MD@dot}\mkern1mu}}%
\def\MD@iddots#1{{\MD@changestyle{#1}%
  \mkern1mu\raise \MD@unit\hbox{\MD@dot}%
  \mkern2mu\raise4\MD@unit\hbox{\MD@dot}%
  \mkern2mu\raise7\MD@unit\vbox{\kern7\MD@unit\hbox{\MD@dot}}}}%
\def\MD@vdots#1{\vbox{\MD@changestyle{#1}%
    \baselineskip4\MD@unit\lineskiplimit\z@
    \kern6\MD@unit\hbox{\MD@dot}\hbox{\MD@dot}\hbox{\MD@dot}}}%
\ifMDlatex
  \DeclareRobustCommand\ddots{\mathinner{\MD@palette\MD@ddots}}%
  \DeclareRobustCommand\iddots{\mathinner{\MD@palette\MD@iddots}}%
  \DeclareRobustCommand\vdots{\mathinner{\MD@palette\MD@vdots}}%
\else
  \def\ddots{\mathinner{\MD@palette\MD@ddots}}%
  \def\iddots{\mathinner{\MD@palette\MD@iddots}}%
  \def\vdots{\mathinner{\MD@palette\MD@vdots}}%
\fi
\catcode`@=12    

\newcommand {\comment}[1]{} 
\newcommand{\C}{{\mathbb C}}
\newcommand{\R}{{\mathbb R}}

\begin{document}
\title{Characterization of the dissipative mappings and their application to perturbations of dissipative-Hamiltonian systems}
\author{Mohit Kumar Baghel\thanks{Department of Mathematics, Indian Institute of Technology Delhi, Hauz Khas, 110016, India. 
Email: \{maz188260, punit.sharma\}@maths.iitd.ac.be. MB acknowledges the support of institute Ph.D. fellowship  by IIT Delhi, India. 
PS acknowledges the support of the DST-Inspire Faculty Award (MI01807-G) by Government of India and FIRP project (FIRP/Proposal Id - 135) by IIT Delhi, India }
\and Nicolas Gillis\thanks{Department of Mathematics and Operational Research, University of Mons, Rue de Houdain 9, 7000 Mons, Belgium. 
Email:  nicolas.gillis@umons.ac.be.
NG acknowledges the support by ERC starting grant No 679515, the Fonds de la Recherche Scientifique - FNRS
and  the  Fonds  Wetenschappelijk  Onderzoek - Vlaanderen (FWO) 
under EOS project O005318F-RG47. } 
\and Punit Sharma$^*$}

\date{}

\maketitle

\begin{abstract}
In this paper, we find necessary and sufficient conditions to identify pairs of matrices $X$ and $Y$ for which there exists $\Delta \in \C^{n,n}$ such that
$\Delta+\Delta^*$ is positive semidefinite and $\Delta X=Y$.
Such a $\Delta$ is called a  
dissipative mapping taking $X$ to $Y$. 
We also provide two different characterizations for the set of all dissipative mappings, and use them to characterize the unique dissipative mapping with minimal Frobenius norm. The minimal-norm dissipative mapping is then used to determine the distance to asymptotic instability for dissipative-Hamiltonian systems under general structure-preserving perturbations. We illustrate our results over some numerical examples and compare them with those of 
Mehl, Mehrmann and Sharma 
(Stability Radii for Linear Hamiltonian 
Systems with Dissipation Under Structure-Preserving Perturbations, SIAM J. Mat. Anal. Appl.\ 37 (4): 1625-1654, 2016). 
\end{abstract}

\section{Introduction}\label{sec:intro}

In this paper, we consider the dissipative mapping problem. 
Given $X,Y \in \C^{n,m}$, the dissipative mapping problem can be divided in three subproblems:  
\begin{enumerate}

\item {\it Existence}: find necessary and sufficient conditions on matrices $X$ and $Y$ for the existence of matrices $\Delta \in \C^{n,n}$ such that $\Delta+\Delta^* $ is positive semidefinite and $\Delta X=Y$.

\item {\it Characterization}: characterize all such dissipative matrices taking $X$ to $Y$. 

\item {\it Minimal norm}: characterize all solutions to the dissipative mapping problem that have minimal norm. In this paper, we focus on the Frobenius norm, which is the standard in the literature.

\end{enumerate} 

The dissipative mapping problem belongs to a wider class of mapping problems: 
for $X,Y \in \C^{n,m}$, the structured mapping problems require to find $\Delta \in \mathcal S \subseteq \C^{n,n}$ such that $\Delta X=Y$, where $\mathcal S$ stands for the structure of the mapping. 
The mapping problems have been extensively studied in~\cite{MacMT08} and~\cite{Adh08,AdhA16} for the structures that are associated with orthosymmetric scalar products. 
These include symmetric, skew-symmetric, Hermitian, skew-Hermitian, Hamiltonian, skew-Hamiltonian, Hermitian, skew-Hermitian, persymmetric, and per-skew symmetric matrices. If $\mathcal S \subseteq \R^{n,n}$,  then it is the real structured mapping problem. In particular, for given $X,Y \in \C^{n,m}$, the real dissipative mapping problem is to find $\Delta \in \R^{n,n}$ such that $\Delta+\Delta^T$ is positive semidefinite and $\Delta X=Y$.

The structured mapping problems occur and are useful in solving various distance problems related to structured matrices and matrix 
polynomials, see for example~\cite{BorKMS14,BorKMS15,MehMS16,MehMS17,OveV05} and the references therein. 
A  mapping problem closely related to the dissipative mapping is the positive semidefinite (PSD) mapping, when $\mathcal S$ is the cone of Hermitian positive semidefinite matrices. 
The PSD mapping problem has been recently solved and used in~\cite{MehMS16,MehMS17}  to derive formulas for structured distances to instability for linear time-invariant dissipative Hamiltonian (DH) systems. 
DH systems have the form 
\begin{equation}\label{eq:defDH1}
\dot{x}(t)=(J-R)Qx(t),
\end{equation}
where $x(t)$ is the state vector at time $t$, 
$J^*=-J \in \mathbb F^{n,n}$, 
$R^*=R \in \mathbb F^{n,n}$ is positive semidefinite and is referred to as the dissipation matrix, and $Q=Q^* \in \mathbb F^{n,n}$ is positive definite and  describes the energy of the system via the function $x \mapsto  x^*Qx$.  When the coefficient matrices in~\eqref{eq:defDH1} are complex, that is,. $\mathbb F=\C$, then it is called a complex DH system, and when the coefficient matrices in~\eqref{eq:defDH1} are real, that is, $\mathbb F=\R$, then it is called a real DH system. 
Such systems are special cases of port-Hamiltonian systems; see for example~\cite{GolSBM03,Sch06,SchM13,MehMS16,MehMS17}.  
A linear time-invariant (LTI) control system $\dot{x}=Ax$ with $A\in \C^{n,n}$ is called  \textit{asymptotically stable} around its equilibrium point at the origin if 
(i)~for any $\epsilon >0$, there exists $\delta_1 > 0$ 
such that if $\|x(t_0)\|< \delta_1$, then $\|x(t)\| < \epsilon$, for all $t > t_0$, and 
(ii) there exists $\delta_2 > 0$ such that if  $\|x(t_0)\|< \delta_2$, then $x(t)\rightarrow 0$ as $t \rightarrow \infty$. It is called \textit{stable} if only (i) holds. An equivalent algebraic characterization of stability is given in terms of spectral conditions on $A$: the system $\dot{x}=Ax$ is called stable if all eigenvalues of the matrix $A$ are in the closed left half of the complex plane and those on the imaginary axis are semisimple. It is  asymptotically stable if all eigenvalues of $A$ are in the open left half of the complex plane.

Any LTI stable system can be represented in the form~\eqref{eq:defDH1} of a DH system~\cite{GilS16}. DH systems are always stable and remain stable as long as the perturbations preserve the DH structure. Still, they may have eigenvalues on the imaginary axis, i.e., they are not necessarily asymptotically stable. Therefore it is useful to know the \emph{distance to asymptotic instability} for DH systems~\cite{MehMS16,MehMS17} defined as the smallest norm of the perturbation that makes system lose its asymptotic stability.
 Studying this question is an essential topic in many applications, like in power system and circuit simulation, see, e.g., \cite{MarL90,MarPJ07,Mar86}, and multi-body systems, see, e.g., \cite{GraMQSW16,Sch90b}.
This is also useful in some specific applications, 
e.g., the  analysis of disk brake squeal \cite[Example 1.1]{MehMS16}, of mass-spring-damper dynamical systems \cite[Example 4.1]{MehMS17}, and of circuit simulation and power system modelling \cite[Example 1.2]{MehMS16}, where the interest is in studying the stability or the instability under perturbation of the matrices 
$J$, $R$ and $Q$. In the DH modelling of physical systems, the matrix $J$ describes the energy flux among the system's energy storage elements, and $R$ represents the damping effects in the system. Thus perturbing $J$
and $R$ together, or only one at a time, is of particular interest. 

In~\cite{MehMS16},  the authors derived various distances to asymptotic instability for DH systems~\eqref{eq:defDH1} while focusing on perturbations  that affect only one matrix from $\{J,R,Q\}$ at a time. 
 Similarly in~\cite{MehMS17}, the authors derived real distances to instability for real DH systems while perturbing only the dissipation matrix $R$.
The framework suggested in~\cite{MehMS16} depends heavily on reformulating the instability radius problem in terms of an equivalent problem of minimizing the generalized Rayleigh quotient of two positive semidefinite matrices.  This  reformulation was achieved using minimal norm skew-Hermitian mappings in case only $J$ is perturbed, 
and PSD mappings in case only $R$ or $Q$ is perturbed. 
However, we note that the framework suggested in~\cite{MehMS16} does not work if we allow perturbations in the DH system that affect more than one matrix from $\{J,R,Q\}$ at a time. 
Analyzing the stability of the system~\eqref{eq:defDH1} when more than one matrix from $\{J,R,Q\}$ is  perturbed  is one of the main motivations  of our work. 

In this paper, we focus on perturbations of DH systems that affect both $J$ and $R$ simultaneously.  More precisely, we find (see Section~\ref{seac:DH}) that the minimal-norm solution to the dissipative mapping problem can be a necessary tool in computing the  structured distance to asymptotic instability for DH systems~\eqref{eq:defDH1} when both $J$ and $R$ are simultaneously perturbed.

\subsection{Contributions and outline of the paper}

In Section 2, we present some preliminary results that will be needed to solve the 
dissipative mapping problem. For the solutions to the
dissipative mapping problem, we present two different characterizations in terms of three matrix variables $K,G,Z$ with symmetry and semidefinite structures. Both  characterizations have an advantage over each other. The first characterization (see Section 3) results in a straightforward computation of the minimal-norm solutions to the dissipative mapping problem, but the  matrix variables  $K,G,Z$ are highly constrained. We also derive necessary and sufficient conditions for solving the real dissipative mapping problem and compute solutions that are of minimal Frobenius norm. 
The second characterization (see Section 4) has the advantage of having a simple form in terms of matrix variables $K,G$, and $Z$. On the other hand, it is unclear how to find minimal-norm solutions via this second characterization. 
The minimal-norm dissipative mapping is  used in Section 5 for studying the structured distance to instability for DH systems~\eqref{eq:defDH1} for simultaneous perturbations of $J$ and $R$.
In Section 6, we present numerical examples comparing these distances to instability  with those of~\cite{MehMS16} where  perturbations affect only one of 
the matrices $J$, $R$, or $Q$. 
 
\paragraph{Notation} In the following, we denote the identity matrix of size $n \times n$ by $I_n$, the spectral norm of a matrix or a vector by $\|\cdot\|$, and the Frobenius norm by ${\|\cdot\|}_F$. The Moore-Penrose pseudoinverse of a matrix or a vector $X$ is denoted by $X^\dagger$, 
and  $\mathcal P_X=I_n-XX^\dagger$ denotes the orthogonal projection onto the null space of $n \times n$ matrix $X^*$. For a square matrix $A$, its Hermitian and skew-Hermitian parts are respectively denoted by $A_H=\frac{A+A^*}{2}$ and $A_S=\frac{A-A^*}{2}$.
For $A=A^* \in \mathbb F^{n,n}$, where $\mathbb F \in \{\R,\C\}$, we denote
$A \succ 0$ ($A\prec 0$) and $A\succeq 0$ ($A\preceq$) if $A $ is Hermitian positive definite (negative definite) and Hermitian positive semidefinite (negative semidefinite), respectively, 
and $\Lambda(A)$ denotes the set of all eigenvalues of the matrix $A$.
For a given matrix $A\in \C^{n,m}$, we use the term {SVD} to denote the standard singular value decomposition of $A$, and {reduced SVD} for the decomposition $A=U_1\Sigma_1V_1^*$ in which $\Sigma_1$ is a square diagonal matrix of size equal to rank of $A$, $r$, with the nonzero singular values of $A$ on its diagonal, and $U_1 \in \C^{n,r}$ and $V_1 \in \C^{m,r}$ have orthonormal columns. 

\section{Preliminaries}

In this section, we discuss some basic results from the literature
and  derive some elementary lemmas that will be necessary to solve the dissipative mapping problem. 

Let us start with two well-known lemmas for positive semidefinite matrices.  

\begin{lemma}{\rm \cite{HorJ85}}\label{psdsimilar}
Let $P\in \C^{n,n}$ and $X \in \C^{n,m}$. If $P\succeq 0$, then $X^*PX \succeq 0$.
\end{lemma}
%

\begin{lemma}{\rm \cite{Alb69}} \label{lem:psd_character}
Let the integer $s$ be such that $0<s<n$, and $R=R^* \in \C^{n,n}$ be partitioned as
$R=\mat{cc}B & C^*\\C & D\rix$ with $B\in \C^{s,s}$, $C\in \C^{n-s,s}$ and $D \in \C^{n-s,n-s}$. Then $R \succeq 0$ if and only if
\begin{enumerate}
\item $B \succeq 0$,
\item $\operatorname{null}(B) \subseteq \operatorname{null}(C)$, and
\item $D-CB^{\dagger}C^* \succeq 0$, where $B^{\dagger}$ denotes the Moore-Penrose pseudoinverse of $B$.
\end{enumerate}
\end{lemma}

Next, we state a result~\cite[Theorem 2.2.3]{Adh08} for skew-Hermitian mappings in terms that allow a direct use in deriving the second characterization for dissipative mappings. 

\begin{theorem}\label{thm:skew-herm-map}
Let $X,Y \in \C^{n,k}$ and define $\mathcal S:=\{\Delta \in \C^{n,n}:~\Delta^*=-\Delta,\,\Delta X=Y\}$. Then $\mathcal S \neq \emptyset$ if and only if $(X^*Y)^*=-X^*Y$. Further, if $\mathcal S \neq \emptyset$, then
\begin{equation*}
\mathcal S = \left\{
YX^\dagger-(YX^\dagger)^*+(X^\dagger)^*X^*YX^\dagger + \mathcal P_XZ\mathcal P_X:~Z\in \C^{n,n},\,Z^*=-Z
\right\}.
\end{equation*}
\end{theorem}

The next two lemmas  will be used in deriving real dissipative matrices taking a complex $X \in \C^{n,m}$ to a complex  $Y \in \C^{n,m}$.

\begin{lemma}{\rm\cite[Lemma 3.3]{AlaBKMM11}}\label{realaux}
Let $A,B \in \C^{n,p}$. Then $[A~\overline A][B~\overline B]^{\dagger}$ is a real matrix.
\end{lemma}

\begin{lemma} \label{real_lemma2}
	Let $X,Y \in \mathbb{C}^{n,k}$ be such that $M=([X ~\overline{X}]^*[Y ~\overline{Y}] + [Y ~\overline{Y}]^*[X ~\overline{X}])^{-1}$  exists. Then $[X ~\overline{X}]M [Y ~\overline{Y}]^*$ is a real matrix.
\end{lemma}
\proof
Let $Q= \begin{bmatrix} O & I_k\\ I_k & O \end{bmatrix}$, where $O \in \C^{k,k}$  is the zero matrix. Then 
\begin{align*}
\overline{[X ~ \overline{X}]M[Y ~ \overline{Y}]^*} = \overline{[X ~ \overline{X}]Q^2 M Q^2[Y ~ \overline{Y}]^*}
= \overline{([X ~ \overline{X}]Q) (Q M Q) (Q [Y ~ \overline{Y}]^*)}=[X ~ \overline{X}] \overline{Q M Q} [Y ~\overline{Y}]^*,
\end{align*}
because $[X ~\overline{X}]Q= [\overline{X}~ X]$ and $Q [Y ~ \overline{Y}]^*= [\overline{Y} ~ Y]^*$. 
Thus we prove the assertion by showing that 
$ \overline{QMQ}=M$. This follows from the following derivations: 
\begin{align*}
\overline{Q M Q} &= \overline{Q([X~ \overline{X}]^*[Y \overline{Y}]+ [Y~ \overline{Y}]^*[X ~ \overline{X}])^{-1} Q } \\
&= Q \overline{([X~ \overline{X}]^*[Y ~\overline{Y}]+ [Y~ \overline{Y}]^*[X ~\overline{X}])^{-1}  } Q\\
&= Q \overline{([X~ \overline{X}]^*[Y ~\overline{Y}]+ [Y~ \overline{Y}]^*[X ~\overline{X}])  }^{-1} Q\\
&= Q ([\overline{X}~ X]^*[ \overline{Y}~ Y]+ [ \overline{Y} ~ Y]^*[ \overline{X} ~ X])^{-1}   Q\\
&= \big (Q [\overline{X}~ X]^*[ \overline{Y}~ Y]Q+ Q[ \overline{Y} ~ Y]^*[ \overline{X} ~ X]   Q\big )^{-1}\\
&= \big ([X ~ \overline{X}]^* [Y ~ \overline{Y}] + [Y ~ \overline{Y}]^*[X \overline{X}]\big)^{-1}\\
&= M.
\end{align*}
\eproof

We close the section by discussing  some simple mapping results that will be necessary to compute the structured distances to instability in Section~\ref{seac:DH}.

\begin{lemma}\label{lem:norms_eq}
Let $x \in \C^{q}$ and $y \in \C^{r}$. Then
\begin{eqnarray}\label{eq:lemt1}
\inf\left\{{\|\Delta_1\|}_F^2+{\|\Delta_2\|}_F^2:~\Delta_1,\Delta_2 \in \C^{r,q}, (\Delta_1-\Delta_2)x=y\right\}=
\inf\left\{\frac{{\|\Delta\|}_F^2}{2}:~\Delta\in \C^{r,q}, \Delta x=y\right\}.
\end{eqnarray}
\end{lemma}
\proof First, let $\Delta \in \C^{r,q}$ such that $\Delta x=y$, and set $\Delta_1=\Delta/2$ and $\Delta_2=-\Delta/2$. Then 
$(\Delta_1-\Delta_2)x=\Delta x=y$,
and ${\|\Delta_1\|}_F^2+{\|\Delta_2\|}_F^2=\frac{{\|\Delta\|}_F^2}{4}+\frac{{\|\Delta\|}_F^2}{4}=\frac{{\|\Delta\|}_F^2}{2}$. This implies ``$\leq$"
in~\eqref{eq:lemt1}.
Conversely, let $\Delta_1,\Delta_2 \in \C^{r,q}$ be such that $(\Delta_1-\Delta_2)x=y$, and set $\Delta=\Delta_1-\Delta_2$. Then 
$\Delta x=y$ and
\begin{eqnarray*}
{\|\Delta \|}_F^2={\|\Delta_1-\Delta_2 \|}_F^2 &\leq& \left({\|\Delta_1\|}_F+{\|\Delta_2\|}_F\right)^2 \\
&=& {\|\Delta_1\|}_F^2+{\|\Delta_2\|}_F^2+2{\|\Delta_1\|}_F{\|\Delta_2\|}_F \\
&\leq&2\left({\|\Delta_1\|}_F^2+{\|\Delta_2\|}_F^2\right).
\end{eqnarray*}
Thus $\frac{{\|\Delta\|}_F^2}{2} \leq {\|\Delta_1\|}_F^2+{\|\Delta_2\|}_F^2 $. This implies ``$\geq$"
in~\eqref{eq:lemt1}.
\eproof

\begin{lemma}\label{lem:norms_eq_ssh}
Let $x,y \in \C^{r}$. Then
\begin{eqnarray}\label{eq:lemt2}
&\inf\left\{{\|\Delta_1\|}_F^2+{\|\Delta_2\|}_F^2:~\Delta_1,\Delta_2\in \C^{r,r},\Delta_1^*=-\Delta_1,\Delta_2 \preceq 0 , (\Delta_1-\Delta_2)x=y\right\}\nonumber \\
&=\inf\left\{{\|\Delta\|}_F^2:~\Delta\in \C^{r,r},\Delta+\Delta^*\succeq 0, \Delta x=y\right\}.
\end{eqnarray}
\end{lemma}
\proof 
The idea behind the proof is similar to Lemma~\ref{lem:norms_eq}. Let $\Delta\in \C^{r,r}$ be such that $\Delta+\Delta^*\succeq 0 $ and $\Delta x=y$. Set $\Delta_1= \Delta_S$ and $\Delta_2= -\Delta_H$. Then clearly $\Delta_1^*=-\Delta_1$, $\Delta_2 \preceq 0$ since $\Delta+\Delta^*\succeq 0$, and $(\Delta_1-\Delta_2)x=\Delta x=y$. Also 
${\|\Delta_1\|}_F^2+{\|\Delta_2\|}_F^2={\left\| \Delta_S\right\|}_F^2+{\left\| \Delta_H \right\|}_F^2={\|\Delta\|}_F^2$.  This implies ``$\leq$"
in~\eqref{eq:lemt2}. Conversely, let $\Delta_1,\Delta_2 \in \C^{r,r}$ be such that
$\Delta_1^*=-\Delta_1$, $\Delta_2 \preceq 0$, and $(\Delta_1-\Delta_2)x=y$. Now set $\Delta=\Delta_1-\Delta_2$. Then clearly $\Delta x=y$ and $\Delta+\Delta^* = -2 \Delta_2 \succeq 0$, since $\Delta_2 \preceq 0$. Further we have that  ${\|\Delta\|}_F^2={\left\| \Delta_S\right\|}_F^2+{\left\| \Delta_H \right\|}_F^2= {\|\Delta_1\|}_F^2+{\|\Delta_2\|}_F^2$. This proves ``$\geq$"
in~\eqref{eq:lemt2}. 
\eproof


\begin{lemma}{\rm \cite[Theorem 2.8]{MehMS16}}\label{lem:auxmap}
Let $B \in \C^{n,r}$ with $\text{{\rm rank}}(B)=r$, $x \in \C^{r}\setminus \{0\}$,   $y \in \C^{n}\setminus \{0\}$, and let $\Delta \in \C^{r,r}$.
Then $B\Delta x=y$ if and only if $\Delta x=B^\dagger y$ and $BB^\dagger y=y$.
\end{lemma}

\section{First characterization of dissipative mappings}\label{sec:mappings}

In this section, we derive the first characterization of dissipative mappings. 
This characterization allows us to find the minimal Frobenius norm solutions to 
the dissipative mapping problem which is shown to be unique. This minimal-norm solution will turn out to be a necessary tool in computing the structured stability radius in Section~\ref{seac:DH}. For given $X \in \C^{n,m}$ and  $Y \in \C^{n,m}$, we define the set of dissipative mappings from $X$ to $Y$ as follows 
\[
\mathbb S(X,Y):=\{\Delta \in \C^{n,n}~:~\Delta+\Delta^* \succeq 0,~\Delta X=Y\}.
\]

We will need the following lemma in characterizing the set of all solutions to the dissipative mapping problem.

\begin{lemma}\label{lem:psdequalcond}
Let $X,Y \in \C^{n,m}$. Suppose that  $\text{{\rm rank}}(X)=p$ and consider the reduced singular value decomposition $X=U_1\Sigma_1 {V}_1^*$ with
$U_1 \in \C^{n,p}$, $\Sigma_1 \in \C^{p,p}$ and $V_1 \in \C^{m,p}$. If $X^*Y+Y^*X \succeq 0$, then $U_1^*\left(YX^\dagger + {(YX^\dagger)}^*\right)U_1 \succeq 0$.
\end{lemma}
\proof Since $X=U_1\Sigma_1 {V}_1^*$, we have $X^\dagger=V_1 \Sigma_1^{-1}U_1^*$. Thus 
\begin{eqnarray*}
X^*Y+Y^*X \succeq 0 
&\iff & V_1\Sigma_1U_1^*Y+Y^*U_1\Sigma_1V_1^* \succeq 0 \\
&\Longrightarrow & V_1^*\left(V_1\Sigma_1U_1^*Y+Y^*U_1\Sigma_1V_1^*\right)V_1 \succeq 0 \quad (\because\,\text{using\,Lemma}\,\ref{psdsimilar})\\
&\iff & \Sigma_1U_1^*YV_1+V_1^*Y^*U_1\Sigma_1 \succeq 0 \hspace{1.5cm}  (\because V_1^*V_1=I_p)\\
&\iff & \Sigma_1^{-1}\left(\Sigma_1U_1^*YV_1+V_1^*Y^*U_1\Sigma_1\right)\Sigma_1^{-1} \succeq 0 \\
&\iff & U_1^*YV_1\Sigma_1^{-1}+\Sigma_1^{-1}V_1^*Y^*U_1 \succeq 0 \\
&\iff & U_1^*\left(YV_1\Sigma_1^{-1}U_1^*+U_1\Sigma_1^{-1}V_1^*Y^*\right)U_1 \succeq 0 \\
&\iff & U_1^*\left(YX^\dagger + {(YX^\dagger)}^*\right)U_1 \succeq 0,
\end{eqnarray*}
which completes the proof.
\eproof

\begin{theorem}\label{thm1:mapping}
Let $X,Y \in \C^{n,m}$, 
and suppose that ${\rm rank}(X)=p$. Then $\mathbb S(X,Y) \neq \emptyset$ if and only if $YX^\dagger X=Y$ and $X^*Y+Y^*X \succeq 0$.
Moreover, if $\mathbb S(X,Y) \neq \emptyset$, 
then 
\begin{enumerate}
\item {\rm Characterization:}~Let  $X=U\Sigma V^*$ be the the singular value decomposition of $X$ with $U=\left[U_1 ~ U_2\right]$, where $U_1 \in \C^{n,p}$. Then 

\begin{eqnarray}\label{eq:S_char}
\mathbb S(X,Y)=\Big\{U\mat{cc}U_1^*YX^\dagger U_1 & U_1^*(YX^\dagger)^*U_2+U_1^*ZU_2\\U_2^*YX^\dagger U_1 & U_2^*(K+G)U_2 \rix U^*:\nonumber \\
Z,K,G \in \C^{n,n}~\text{satisfying}~\eqref{eq:cond1}-\eqref{eq:cond3} \Big\},
\end{eqnarray} 
%
where
\begin{eqnarray}
&&G^*=-G,\quad K\succeq 0, \label{eq:cond1}\\
&&{\rm null}\left(U_1^*(YX^\dagger+(YX^\dagger)^*)U_1\right) \subseteq {\rm null}\left(U_2^*(2YX^\dagger +Z^*)U_1\right), \label{eq:cond2}\\
&&K-\frac{1}{8}(2YX^\dagger+Z^*){\left(XX^\dagger YX^\dagger+(YX^\dagger)^*XX^\dagger\right)}^{\dagger}(2YX^\dagger+Z^*)^* \succeq 0\label{eq:cond3}.
\end{eqnarray}

\item {\rm Minimal norm mapping:}~
\begin{equation}\label{eq:minres}
\inf_{\Delta \in \mathbb S(X,Y)} {\|\Delta\|}_F^2 \; = \; 2{\|YX^\dagger\|}_F^2-{\rm trace}\left((YX^\dagger)^*XX^\dagger(YX^\dagger)\right),
\end{equation}
where the infimum is uniquely attained by the matrix $\mathcal H:=YX^\dagger -{(YX^\dagger)}^* \mathcal P_X$,
which is obtained by setting $K=0$, $G=0$, and $Z=-2U_1^*(YX^\dagger)^*U_2$ in~\eqref{eq:S_char}.
\end{enumerate} 
\end{theorem}
\proof First suppose that  $\Delta \in \mathbb S(X,Y)$, i.e., $\Delta+\Delta^* \succeq 0$ and $\Delta X=Y$. Then 
$YX^\dagger X= \Delta X X^\dagger  X= \Delta X=Y$. 
By Lemma~\ref{psdsimilar}
\[
X^*Y+Y^*X=X^*\Delta X+X^* \Delta^* X=X^*(\Delta +\Delta^*)X \succeq 0,
\]
since $\Delta+\Delta^* \succeq 0$. Conversely, suppose that  $X$ and $Y$ satisfy $YX^\dagger X=Y$ and $X^*Y+Y^*X \succeq 0$. Then the matrix
$\mathcal H$ satisfies $\mathcal H X=YX^\dagger X=Y$ and
\begin{eqnarray*}
\mathcal H +{\mathcal H}^*&=&YX^\dagger -{(YX^\dagger)}^* \mathcal P_X + {(YX^\dagger)}^* -\mathcal P_X {YX^\dagger}\\
&=& {(YX^\dagger)}^*XX^\dagger + XX^\dagger (YX^\dagger)\quad \quad (\because \mathcal P_X=I_n-XX^\dagger)\\
&=& XX^\dagger \left({(YX^\dagger)}^*+YX^\dagger\right)XX^\dagger \quad \quad (\because X^\dagger XX^\dagger=X^\dagger,~(XX^\dagger)^*=XX^\dagger)\\
&=& U_1U_1^* \big({(YX^\dagger)}^*+YX^\dagger\big)U_1U_1^* \quad \quad (\because  XX^\dagger=U_1U_1^*)\\
&=& \mat{cc}U_1&U_2\rix \mat{cc}U_1^* \big({(YX^\dagger)}^*+YX^\dagger\big)U_1 &0 \\0& 0\rix \mat{cc}U_1&U_2 \rix^*\\
&\succeq & 0,
\end{eqnarray*}
where the last identity follows by using Lemma~\ref{lem:psdequalcond}, since $X^*Y+Y^*X \succeq 0$, and by Lemma~\ref{psdsimilar}.

Next, we prove~\eqref{eq:S_char}. First suppose that $\Delta \in \mathbb S(X,Y)$, i.e., $\Delta+\Delta^* \succeq 0$ and $\Delta X=Y$.
Let $\Sigma=\mat{cc}\Sigma_1 & 0\\ 0 &0\rix$ and $V=\mat{cc}V_1 & V_2 \rix$, where $V_1 \in \C^{m,p}$, $V_2 \in \C^{m, m-p}$, and $\Sigma_1 \in \C^{p,p}$
such that $X=U_1\Sigma_1V_1^*$ is the reduced SVD of $X$. Then $X^*=V_1\Sigma_1U_1^*$ and $X^\dagger=V_1 \Sigma_1^{-1}U_1^*$.
Consider $\Delta=UU^*\Delta UU^*$ and $\widetilde \Delta=U^*\Delta U=\widetilde \Delta_H+\widetilde \Delta_S$, where
\[
\widetilde \Delta_H=U^*\Delta_H U=\mat{cc}H_{11}&H_{12}\\H_{12}^* & H_{22}\rix\quad \text{and}\quad
\widetilde \Delta_S=U^*\Delta_S U=\mat{cc}S_{11}&S_{12}\\-S_{12}^* & S_{22}\rix.
\]
Clearly, ${\|\Delta\|}_{F}={\|\widetilde \Delta\|}_{F}$ and also $\Delta_H \succeq 0$ $\Longleftrightarrow$ $\widetilde \Delta_H \succeq 0$.
As $\Delta X=Y$, we have
\begin{eqnarray*}
U^*\Delta U U^* X=U^*Y&\Longrightarrow& \widetilde \Delta \mat{c}U_1^* \\U_2^* \rix X= \mat{c}U_1^*Y \\U_2^*Y \rix \\
 &\Longrightarrow&
\mat{cc}H_{11}+S_{11} & H_{12}+S_{12}\\ H_{12}^*-S_{12}^* & H_{22}+S_{22}\rix \mat{c}\Sigma_1V_1^* \\ 0\rix=\mat{c}U_1^*Y \\U_2^*Y \rix.
\end{eqnarray*}
This implies that
\begin{equation}\label{eq:h11}
\left(H_{11}+S_{11} \right)\Sigma_1V_1^*=U_1^*Y,
\end{equation}
and
\begin{equation}\label{eq:h12}
\left(H_{12}^*-S_{12}^* \right)\Sigma_1V_1^*=U_2^*Y.
\end{equation}
Thus from~\eqref{eq:h11}, we have $H_{11}+S_{11}=U_1^*YV_1\Sigma_1^{-1}=U_1^*YX^\dagger U_1$,
since $X^\dagger=V_1\Sigma_1^{-1}U_1^*$ and $X^\dagger U_1=V_1 \Sigma_1^{-1}$. This implies that
\begin{equation}\label{eq:h11form}
H_{11}=U_1^*\left(\frac{(YX^\dagger)+(YX^\dagger)^*}{2}\right)U_1 \quad \text{and} \quad
S_{11}=U_1^*\left(\frac{(YX^\dagger)-(YX^\dagger)^*}{2}\right)U_1.
\end{equation}
Note that since $X^*Y+Y^*X \succeq 0$, in view of Lemma~\ref{lem:psdequalcond}, we have that $H_{11} \succeq 0$.
Similarly, from~\eqref{eq:h12} we have $H_{12}^*-S_{12}^*=U_2^*YV_1\Sigma_1^{-1}=U_2^*YX^\dagger U_1$. This implies that
\begin{equation}\label{eqh12form}
H_{12}=U_1^*(YX^\dagger)^* U_2 + S_{12},
\end{equation}
where $S_{12} \in \C^{p,n-p}$ is a matrix variable. Thus from~\eqref{eq:h11form} and~\eqref{eqh12form}, $\widetilde \Delta$ has the form
\begin{equation}\label{eq:delndelt}
\widetilde \Delta= \mat{cc}H_{11}+S_{11} & H_{12}+S_{12}\\ H_{12}^*-S_{12}^* & H_{22}+S_{22}\rix=
\mat{cc} U_1^*YX^\dagger U_1 & U_1^*(YX^\dagger)^*U_2 + 2S_{12} \\ U_2^*YX^\dagger U_1 & H_{22}+S_{22}\rix,
\end{equation}
where $H_{22},S_{22} \in \C^{n-p,n-p}$ and $S_{12} \in \C^{p,n-p}$ are such that $\widetilde \Delta_H \succeq 0$ and $\widetilde \Delta_S=-\widetilde \Delta_S^*$.
That means, in view of Lemma~\ref{lem:psd_character}, $S_{22}$ satisfies that $S_{22}^*=-S_{22}$, and $H_{22}$ and $S_{12}$
satisfy the following constraints: $H_{22}\succeq 0$,
\begin{equation}\label{eq:H22}
H_{22}-(U_2^*YX^\dagger U_1+S_{12}^*){\left(\frac{U_1^*(YX^\dagger+(YX^\dagger)^*)U_1}{2}\right)}^{\dagger}(U_2^*YX^\dagger U_1+S_{12}^*)^* \succeq 0,
\end{equation}
and
\begin{equation}\label{eq:nullS12}
\text{null}\left(\frac{U_1^*(YX^\dagger+(YX^\dagger)^*)U_1}{2}\right) \subseteq \text{null}\left(U_2^*YX^\dagger U_1+S_{12}^*\right).
\end{equation}
Thus from~\eqref{eq:delndelt}, we have

\begin{eqnarray}\label{eq:const}
\Delta&=&U\widetilde \Delta U^*=\mat{cc}U_1 & U_2\rix \mat{cc} U_1^*YX^\dagger U_1 & U_1^*(YX^\dagger)^*U_2 + 2S_{12} \\ U_2^*YX^\dagger U_1 & H_{22}+S_{22}\rix
\mat{c}U_1^* \\ U_2^*\rix.
\end{eqnarray}
By setting $Z=2U_1S_{12}U_2^*$, $K=U_2H_{22}U_2^*$, $G=U_2S_{22}U_2^*$, and by using the fact that $U_1U_1^*+ U_2U_2^*=UU^*=I_n$, $U_1U_1^*=XX^\dagger$, $U_1^*U_1=I_p$, and $U_2^*U_2=I_{n-p}$, we obtain that 
\begin{equation*}
\Delta=U\mat{cc}U_1^*YX^\dagger U_1 & U_1^*(YX^\dagger)^*U_2+U_1^*ZU_2\\U_2^*YX^\dagger U_1 & U_2^*(K+G)U_2 \rix U^*,
\end{equation*}
where $G$, $Z$, and $K$ satisfy the conditions~\eqref{eq:cond1}--\eqref{eq:cond3}. This proves $``\subseteq"$ in~\eqref{eq:S_char}.

For the other  inclusion in~\eqref{eq:S_char}, let

\begin{equation*}
A=U\mat{cc}U_1^*YX^\dagger U_1 & U_1^*(YX^\dagger)^*U_2+U_1^*ZU_2\\U_2^*YX^\dagger U_1 & U_2^*(K+G)U_2 \rix U^*,
\end{equation*}
where $G$, $Z$, and $K$ satisfy the conditions~\eqref{eq:cond1}--\eqref{eq:cond3}, which can be written as 
\begin{eqnarray}
A=YX^\dagger+(YX^\dagger)^* \mathcal P_X+XX^\dagger Z\mathcal P_X+\mathcal P_X K \mathcal P_X + \mathcal P_X G \mathcal P_X.
\end{eqnarray}
Clearly  $AX=Y$ since $YX^\dagger X=Y$ and $\mathcal P_XX=0$. Also $A+A^* \succeq 0$. Indeed,
\begin{equation*}
U^*AU=
\mat{cc}
U_1^*YX^\dagger U_1 & U_1^*(YX^\dagger)^*U_2+U_1^*ZU_2\\U_2^*YX^\dagger U_1 & U_2^*(K+G)U_2 \rix,
\end{equation*}
and thus
\[
U^*(A+A^*)U=\mat{cc}
U_1^*\left(YX^\dagger+(YX^\dagger)^*\right) U_1 & 2U_1^*(YX^\dagger)^*U_2 +U_1^*ZU_2\\
2U_2^*YX^\dagger U_1 +U_2^*Z^*U_1 & 2U_2^*KU_2 
\rix \succeq 0,
\]
because of Lemma~\ref{lem:psd_character}, since $K$, $Z$, $G$ satisfy~\eqref{eq:cond1}--\eqref{eq:cond3}. This completes the proof of~\eqref{eq:S_char}.

Suppose that $\mathbb S(X,Y) \neq  \emptyset$ and  let $\Delta \in \mathbb S(X,Y)$, then from~\eqref{eq:delndelt} $\Delta$ satisfies
\[
{\|\Delta\|}_F^2 = {\|\widetilde \Delta \|}_F^2 =
{\left\| \mat{cc} U_1^*YX^\dagger U_1 & U_1^*(YX^\dagger)^*U_2 + 2S_{12} \\ U_2^*YX^\dagger U_1 & H_{22}+S_{22}\rix\right\|}_F^2,
\]
where $S_{22}\in \C^{n-p,n-p}$ satisfies that $S_{22}^*=-S_{22}$, and $H_{22}\in \C^{n-p,n-p}$ and $S_{12}\in \C^{p,n-p}$
 satisfy~\eqref{eq:H22} and~\eqref{eq:nullS12}. This implies that
 \begin{eqnarray}
 {\|\Delta\|}_F^2 &=&
 {\|U_1^* YX^\dagger U_1\|}_F^2 +  {\|U_2^* YX^\dagger U_1\|}_F^2 +
  {\|U_1^* (YX^\dagger)^* U_2+2S_{12}\|}_F^2+  {\|H_{22}+S_{22}\|}_F^2 \nonumber\\
  &=& {\| YX^\dagger U_1\|}_F^2 +  {\|U_1^* (YX^\dagger)^* U_2+2S_{12}\|}_F^2+
  {\|H_{22}\|}_F^2+  {\|S_{22}\|}_F^2 ,\label{eq:minn}
 \end{eqnarray}
where the last equality follows as ${\|U_1^* YX^\dagger U_1\|}_F^2 +  {\|U_2^* YX^\dagger U_1\|}_F^2= {\|UYX^\dagger U_1\|}_F^2 = {\|YX^\dagger U_1\|}_F^2$
because $\|\cdot\|_F$ is unitary invariant, and by the fact that for any square matrix $A=A_H+A_S$ we have ${\|A\|}_F^2={\|A_H\|}_F^2+{\|A_S\|}_F^2$.
Further, from~\eqref{eq:minn} we have that
\begin{equation*}
{\|\Delta\|}_F^2 \geq
{\| YX^\dagger U_1\|}_F^2 +  {\|U_1^* (YX^\dagger)^* U_2+2S_{12}\|}_F^2,
\end{equation*}
where the lower bound is attained by setting $H_{22}=0$ and $S_{22}=0$.
Thus any $\Delta \in \mathbb S(X,Y)$ satisfies
\begin{equation*}
{\|\Delta\|}_F^2 \geq
{\| YX^\dagger U_1\|}_F^2 +  {\|U_1^* (YX^\dagger)^* U_2+2S_{12}\|}_F^2,
\end{equation*}
where $S_{12}$ satisfies~\eqref{eq:H22} and~\eqref{eq:nullS12} with $H_{22}=0$, but the only $S_{12}$ that satisfies~\eqref{eq:H22} with $H_{22}=0$
is $S_{12}=-U_1^*(YX^\dagger)^*U_2$ because
$U_1^*\left((YX^\dagger)+(YX^\dagger)^*\right)U_1 \succeq 0$.
Hence by setting $H_{22}=0$, $S_{22}=0$, and $S_{12}=-U_1^*(YX^\dagger)^*U_2$ in~\eqref{eq:delndelt}, we obtain the unique matrix
\begin{eqnarray*}
\Delta&=& \mat{cc}U_1 & U_2\rix
\mat{cc} U_1^*YX^\dagger U_1 & -U_1^*(YX^\dagger)^*U_2  \\ U_2^*YX^\dagger U_1 & 0 \rix
\mat{c} U_1^* \\ U_2^*\rix \\
&=& (U_1U_1^*+U_2U_2^*)YX^\dagger U_1U_1^*
-U_1U_1^*(YX^\dagger)^*U_2U_2^* \\
&=& YX^\dagger XX^\dagger -XX^\dagger(YX^\dagger)^*\mathcal P_X\\
&=&YX^\dagger -(YX^\dagger)^*\mathcal P_X\\
&=& \mathcal H
\end{eqnarray*}
which implies that~\eqref{eq:minres} holds.
\eproof

The following corollary of Theorem~\ref{thm1:mapping} for the vector case ($m=1$) of the minimal norm dissipative mapping  is of particular interest for us as it will be used in Section~\ref{seac:DH} in finding the structured stability radius
for DH systems when both $J$ and $R$ are perturbed. 
\begin{corollary}\label{eq:vectorcase}
Let $x,y \in \C^n \setminus\{0\}$.
Then $\mathbb S(x,y) \neq \emptyset$ if and only if $\real{(x^*y)} \geq 0$. Moreover, if $\mathbb S(x,y)\neq \emptyset$, then
\begin{equation}\label{eq:minres1}
\inf_{\Delta \in \mathbb S(x,y)} {\|\Delta\|}_F^2 \; = \; 2\frac{\|y\|^2}{\|x\|^2}-\frac{|x^*y|^2}{\|x\|^4},
\end{equation}
where the infimum is attained by the unique matrix $\mathcal H:=\frac{yx^*}{\|x\|^2}-\frac{xy^*}{\|x\|^2}+(y^*x)\frac{xx^*}{\|x\|^4}$.
\end{corollary}

\subsection{Real dissipative mappings}\label{sec:real1}
For $X$, $Y \in \C^{n,m}$, if we consider the real dissipative mapping problem, i.e., finding $\Delta \in \R^{n,n}$ such that $\Delta+\Delta^T \succeq 0$ 
and $\Delta X=Y$, then the minimal Frobenius norm solution can be easily obtained from Theorem~\ref{thm1:mapping}. 
To see this, observe that for a real $\Delta$, $\Delta X=Y$ if and only if $\Delta [X~\overline X]=[Y~\overline Y]$. 
In the following, we show that if there exists a complex dissipative mapping $\Delta$ satisfying $\Delta \mathcal X=\mathcal Y$ where $\mathcal X=[X~\overline X]$ and  $\mathcal Y=[Y~\overline Y]$, then there also exists a real dissipative mapping. Moreover, from Theorem~\ref{thm1:mapping}, 
we can easily find a minimal norm real dissipative mapping taking $X$ to $Y$.


\begin{theorem}\label{cor:realthm1:mapping}
Let $X,Y \in \C^{n,m}$, and define  $\mathbb S_{\R}(X,Y):=\{\Delta \in \R^{n,n}~:~\Delta+\Delta^T \succeq 0,~\Delta X=Y\}$. Let $\mathcal X =[X~\overline X]$
and $\mathcal Y=[Y~\overline Y]$.
Suppose that ${\rm rank}(\mathcal X)=p$. Then $\mathbb S_\R(X,Y) \neq \emptyset$ if and only if $\mathcal Y \mathcal X^\dagger \mathcal X=\mathcal Y$ and $\mathcal X^*\mathcal Y+\mathcal Y^*\mathcal X \succeq 0$. Moreover,
if  $\mathbb S_\R(X,Y) \neq \emptyset$ and if we consider the singular value decomposition $\mathcal X=U\Sigma V^T$ with $U=\left[U_1 ~ U_2\right]$, where $U_1 \in \C^{n,p}$, then the set

\begin{equation*}
\left\{ U\mat{cc}U_1^*\mathcal Y \mathcal X^\dagger U_1 & U_1^*(\mathcal Y \mathcal X^\dagger)^*U_2+U_1^*ZU_2\\U_2^*\mathcal Y \mathcal X^\dagger U_1 & U_2^*(K+G)U_2 \rix U^*:Z,K,G \in \R^{n,n}~\text{satisfying}~\eqref{realeq:cond1}-\eqref{realeq:cond3} \right\},
\end{equation*}
where
\begin{eqnarray}
&&G^T=-G,\quad K\succeq 0, \label{realeq:cond1}\\
&&{\rm null}\left(U_1^*(\mathcal Y \mathcal X^\dagger+(\mathcal Y \mathcal X^\dagger)^T)U_1\right) \subseteq {\rm null}\left(U_2^*(2\mathcal Y \mathcal X^\dagger +Z^T)U_1\right), \label{realeq:cond2}\\
&&K-\frac{1}{8}(2 \mathcal Y \mathcal X^\dagger+Z^T){\left(\mathcal X \mathcal X^\dagger \mathcal Y \mathcal X^\dagger+(\mathcal Y \mathcal X^\dagger)^T\mathcal X \mathcal X^\dagger\right)}^{\dagger}(2 \mathcal Y \mathcal X^\dagger+Z^T)^T \succeq 0\label{realeq:cond3},
\end{eqnarray}
 is contained in $\mathbb S_\R(X,Y)$. Further,
\begin{equation*}\label{realeq:minres}
\inf_{\Delta \in \mathbb S_\R(X,Y)} {\|\Delta\|}_F^2 \; = \; 2{\|\mathcal Y \mathcal X^\dagger\|}_F^2-{\rm trace}\left((\mathcal Y
 \mathcal X^\dagger)^*\mathcal X \mathcal X^\dagger(\mathcal Y \mathcal X^\dagger)\right),
\end{equation*}
where the infimum is attained by the matrix $\mathcal H:=\mathcal Y \mathcal X^\dagger -{(\mathcal Y \mathcal X^\dagger)}^* \mathcal P_\mathcal X$.
\end{theorem}
\proof In view of Lemma~\ref{realaux}, the proof follows from Theorem~\ref{thm1:mapping}. Indeed, from Theorem~\ref{thm1:mapping} there exists $\Delta \in \mathbb S(\mathcal X,\mathcal Y)$ if and only if $\mathcal Y \mathcal X^\dagger \mathcal X=\mathcal Y$ and $\mathcal X^*\mathcal Y+\mathcal Y^*\mathcal X \succeq 0$. Assuming the latter two conditions hold true, and consider a family of mappings defined by 
\begin{eqnarray*}
\Delta(Z,K,G):=U\mat{cc}U_1^*\mathcal Y \mathcal X^\dagger U_1 & U_1^*(\mathcal Y \mathcal X^\dagger)^*U_2+U_1^*ZU_2\\U_2^*\mathcal Y \mathcal X^\dagger U_1 & U_2^*(K+G)U_2 \rix U^*,
\end{eqnarray*}
where $Z,K,G \in \R^{n,n}$ satisfy~\eqref{realeq:cond1}-\eqref{realeq:cond3}. Using Lemma~\ref{realaux}, it is easy to check that $\Delta(Z,K,G)$ is real, which implies that $\Delta(Z,K,G) \in \mathbb S_\R( X, Y)$ since $\mathbb S_\R(\mathcal X,\mathcal Y)=\mathbb S_\R(X,Y)$. Further, we have that 
\begin{eqnarray}\label{eq:minmapreal}
\inf_{\Delta \in \mathbb S(\mathcal X,\mathcal Y)} {\|\Delta\|}_F \leq \inf_{\Delta \in \mathbb S_\R(X,Y)} {\|\Delta\|}_F,
\end{eqnarray}
and the left hand side infimum in~\eqref{eq:minmapreal} is attained by the unique map $\mathcal H=\mathcal Y \mathcal X^\dagger -{(\mathcal Y \mathcal X^\dagger)}^* \mathcal P_\mathcal X$. Observe from Lemma~\ref{realaux} that the matrix $\mathcal H$ is real, that is, $\mathcal H \in  \mathbb S_\R(X,Y)$. This implies that 
\begin{equation}
\inf_{\Delta \in \mathbb S_\R(X,Y)} {\|\Delta\|}_F^2= {\|\mathcal H\|}_F^2=2{\|\mathcal Y \mathcal X^\dagger\|}_F^2-{\rm trace}\left((\mathcal Y \mathcal X^\dagger)^*\mathcal X \mathcal X^\dagger(\mathcal Y \mathcal X^\dagger)\right),
\end{equation}
which completes the proof.
\eproof

\section{Second characterization of dissipative mappings}

In this section, for given $X,Y \in \mathbb C^{n,m} \setminus{ \{0\} }$ such that $X^*Y+Y^*X$ is of full rank, we provide a different characterization of dissipative mappings from $X$ to $Y$. The main advantage of this characterization over Theorem~\ref{thm1:mapping}  is that it is explicitly in terms of the matrix variables $K,Z,G \in \C^{n,n}$ such that $K\succeq 0$ and $G^*=-G$, and the matrices $K,Z,G$ do not have to satisfy other constraints like~\eqref{eq:cond1}-\eqref{eq:cond3}. However, extracting the minimal Frobenius norm solutions from this characterization cannot be obtained easily.  
%
\begin{theorem}\label{thm:diffcharac}
Let $X, Y \in \mathbb{C}^{n,m}\setminus{\{0\}}$ and suppose that  $M:=(X^*Y+Y^*X)^{-1}$ exists. 
Then $\mathbb{S}(X,Y) \neq \emptyset $ if and only if $X^*Y+Y^*X \succ 0$ and 
$YX^\dagger X=Y$. Moreover,
\begin{eqnarray} \label{starmatrix}
\mathbb{S}(X,Y) = \{ H + \widetilde{H}(K,G , Z) \mid  K,G, Z\in \mathbb{C}^{n,n},K\succeq 0,~G^*=-G\},
\end{eqnarray}
where $H$ and $\widetilde H$ are respectively defined by
\begin{eqnarray}\label{eqHB}
H :=\frac{1}{2} \left( YX^{\dagger} - X^{\dagger ^*}Y^*  \right)
+  \frac{1}{2}\Big\{ YMY^* + YMX^*YX^{\dagger} 
+ (YX^{\dagger})^* XMY^* \nonumber \\
+ (YX^{\dagger})^* XMX^*YX^{\dagger} \Big\}
\end{eqnarray} 
and
\begin{eqnarray}\label{eqHtilB}
\widetilde{H}(K,G,Z):= \frac{1}{2}\Big({\mathcal P}_X K {\mathcal P}_X + {\mathcal P}_X G {\mathcal P}_X - {\mathcal P}_X Z^*XX^\dagger + X^{\dagger ^ *} X^\dagger Z {\mathcal P}_X +  YMX^* Z {\mathcal P}_X \nonumber \\+ (Y X^\dagger)^*X M X^* Z {\mathcal P}_X  + {\mathcal P}_XZ^* XMY^* 
+  {\mathcal P}_X Z^* XMX^*YX^\dagger+  {\mathcal P}_XZ^*X M X^*Z {\mathcal P}_X\Big).
\end{eqnarray}
\end{theorem}
\proof
If there exists $\Delta$ such that $\Delta X = Y$ and $\Delta + \Delta^* \succeq 0$, 
then
\[
X^*Y+Y^*X=X^*(Y+\Delta^*X)=X^*(\Delta+\Delta^*)X.
\]
This implies that $X^*Y+Y^*X$ is positive definite as $\Delta+\Delta^* \succeq 0$ and $X^*Y+Y^*X$ is invertible. 
For the converse, let $X,Y$ be such that $X^*Y + Y^* X \succ 0$. Then the matrix $H$ in~\eqref{eqHB} is well defined and satisfies 
$HX=Y$. Also $H+H^*$ is a positive semidefinite matrix of the form $BB^*$ for some $B\in \C^{n,n}$.
Indeed, if we let $M=M_1M_1^*$, where $M_1$ is the Cholesky factor of $M$, then we have 
\begin{eqnarray*}
H+H^*&=&YMY^* + YMX^*YX^{\dagger} + (YX^{\dagger})^*XMY^* + (YX^{\dagger})^* XMX^*YX^{\dagger} \\
&=& \big(YM_1 + (YX^\dagger )^*XM_1\big)
\big(YM_1 + (YX^\dagger )^*XM_1\big)^*\\
&=& BB^*,
\end{eqnarray*}
 where  $B=YM_1 + (YX^\dagger )^*XM_1$.

Next we prove $``\subseteq"$ in~\eqref{starmatrix}. For this, let $\Delta \in \mathbb{S}$, 
that is,
$\Delta X=Y$ and $\Delta + \Delta ^* \succeq 0$. 
By~\cite[Lemma 1.3]{Sun93}, there exists $Z\in \mathbb{C}^{n,n}$ such that
\begin{equation}\label{eqdelta}
\Delta = YX^\dagger + Z {\mathcal P}_X.
\end{equation}
Also, we can write $\Delta$ as 
\begin{equation}\label{eqdeltahs}
\Delta= \frac{\Delta + \Delta^*}{2} + \frac{\Delta - \Delta^*}{2}.
\end{equation}
Now since  $\Delta + \Delta^* \succeq 0$ we have $\Delta+ \Delta^* = A^*A $ 
for some $A \in \mathbb{C}^{n,n}$. This implies that 
$(\Delta + \Delta^*)X=A^*AX$ and  thus $\Delta X + \Delta^* X= A^*A X$.
 By setting $W:=AX$, we obtain $Y + \Delta^*X=A^*W$ because $\Delta X=Y.$ This yields the mappings $AX=W$ and $A^*W=Y + \Delta^*X.$
The matrices $X$, $Y$ and $W$ satisfy
\begin{equation}
X^*(Y+\Delta^*X) =X^*Y+X^*\Delta^*X = X^*(\Delta+\Delta^*)X,
=X^*A^*AX=W^*W,
\end{equation}
$YX^{\dagger} X =Y$ and $WX^{\dagger} X = W$. Therefore, from~\cite[Theorem 2.1]{MehMS17}, $A$ can be written as
\begin{equation}
A= WX^{\dagger} +
\left( \left( Y+\Delta^*X \right)W^{\dagger}  \right)^*-
\left(  (Y+\Delta^*X)W^{\dagger}\right)^*XX^\dagger + 
\big( I_n-W W^{\dagger}\big)R_1\mathcal P_X,
\end{equation}
for some $R_1 \in \mathbb{C}^{n,n}$. This can be further simplified as
\begin{equation}\label{simpA}
 A=( ( Y+\Delta^*X )W^{\dagger}  )^* + ( I_n-W W^{\dagger} )R_1 \mathcal P_X,
 \end{equation}
 using the fact that $( (Y+\Delta^*X)W^{\dagger}) )^*XX^\dagger= W X^\dagger$, since $(Y+\Delta^*X )^* X=W^{*} W$.
Thus from~\eqref{simpA}, we have
\begin{eqnarray}\label{tempp1}
A^*A &=& \left( \big( ( Y+\Delta^*X )W ^\dagger  \big)^* + 
\big( I_n-W W^{\dagger} \big)R_1 \mathcal P_X \big)^* \big(\big( ( Y+\Delta^*X )W^{\dagger}  \big)^* + 
\big( I_n-W W^{\dagger} \big)R_1\mathcal P_X \right)  \nonumber \\
&=& (Y+\Delta^* X )M(Y+\Delta^* X)^* + \mathcal P_X R_1^*\big(  I_n-W W^{\dagger}  \big)\big( I_n-W W^{\dagger} \big)R_1\mathcal P_X.
\end{eqnarray}
Define  $K:= R_1^*\big( I_n-W W^{\dagger}\big) \big(  I_n-W W^{\dagger}  \big) R_1$, then  $K \succeq 0$, and from~\eqref{tempp1} we have
\begin{eqnarray}\label{defAstarA}
A^*A= (Y+\Delta^* X ) M (Y+\Delta^* X)^* + \mathcal P_X K \mathcal P_X.
\end{eqnarray}
Now we consider the skew-Hermitian part $\Delta - \Delta^*$ of $\Delta$. Let $S=\Delta - \Delta^*$. Then  $SX=\Delta X- \Delta^* X= Y- \Delta^* X $. Again since
 \[
 X^*(Y-\Delta^*X)=X^*\Delta X-X^*\Delta^* X=X^*(\Delta-\Delta^*)X=
  -(Y^*-X^*\Delta )X= -(Y-\Delta^*X)^*X,
 \]
by Theorem~\ref{thm:skew-herm-map},  we have
\begin{eqnarray}\label{defS}
S=(Y-\Delta^*X)X^{\dagger}- ((Y-\Delta^*X)X^{\dagger})^*\mathcal P_X + \mathcal P_XG\mathcal P_X, 
\end{eqnarray}
for some skew-Hermitian matrix $G \in \mathbb{C}^{n,n}$. By inserting $\Delta=Y X^\dagger + Z \mathcal P_X $ from~\eqref{eqdelta} and $W^{\dagger}(W^{\dagger})^*=((Y+ \Delta ^* X)^* X)^{-1}= (X^*Y + Y^*X)^{-1}$ in~\eqref{defAstarA} and \eqref{defS}, we respectively obtain
\begin{eqnarray}\label{finAstarA}
&A^*A=YMY^* + Y M X^*Y X^{\dagger} + Y M X^* Z \mathcal P_X + (YX^{\dagger})^* X M Y^* + (YX^{\dagger})^*X M X^*YX^{\dagger} \nonumber \\ 
&+ (YX^{\dagger})^* X M X^* Z \mathcal P_X + \mathcal P_X Z^* X M Y^* + \mathcal P_X Z^* X M X^* Y X^{\dagger}+ \mathcal P_X Z^* X M X^* Z \mathcal P_X 
+\mathcal P_XK\mathcal P_X, \nonumber \\
\end{eqnarray}
and
\begin{eqnarray}\label{finS}
S=YX^{\dagger}- X^{\dagger ^ *}Y^*XX^{\dagger} - \mathcal P_X Z^* X X^{\dagger}- X^{\dagger ^ *}Y^*\mathcal P_X + X^{\dagger ^ *}X^{\dagger} Z \mathcal P_X  + \mathcal P_XG\mathcal P_X.
\end{eqnarray}
By inserting $\Delta + \Delta^*= A^*A$ from~\eqref{finAstarA} and $\Delta - \Delta^*=S$ from~\eqref{finS} in~\eqref{eqdeltahs}, and by separating the terms with and without matrices $G$, $K$, $Z$, the matrix $\Delta$ can be written as 
\[ 
\Delta=H+ \widetilde{H}(K,G,Z).
\]
This proves $``\subseteq"$ in~\eqref{starmatrix}.
%

Now, let us prove $``\supseteq"$ in~\eqref{starmatrix}.
Suppose, $\Delta= H+\widetilde{H}(K,G,Z),$ where $H$ is defined 
in~\eqref{eqHB} and $\widetilde{H}(K,G,Z)$ is defined in~\eqref{eqHtilB} for some matrices $G,K,Z\in \C^{n,n}$ such that $K\succeq 0$ and $G^*=-G$. Then 
it is easy to check that $\Delta X=Y$ since $H$ and $\widetilde H$ satisfy
$H X= Y$ and $\widetilde{H}(K,G,Z) X= 0$. Also
\begin{eqnarray}
\Delta + \Delta^*=\big(H+\widetilde{H}(K,G,Z)\big) + \big(H+\widetilde{H}(K,G,Z)\big)^* \succeq 0. 
\end{eqnarray}
Indeed,
\begin{eqnarray}\label{finlast}
&\big(H+\widetilde{H}(K,G,Z)\big) +\big(H+\widetilde{H}(K,G,Z)\big)^*  \nonumber \\
&=(H+H^*)+(\widetilde{H}(K,G,Z)+\widetilde{H}(K,G,Z)^*)\nonumber \\
& =YMY^* + YMX^*YX^{\dagger} + (YX^{\dagger})^*XMY^* + (YX^{\dagger})^* 
  XMX^*YX^{\dagger}   
+ Y M X^* Z \mathcal P_X  +\nonumber  \\
&(YX^{\dagger})^* X M X^* Z \mathcal P_X + \mathcal P_X Z^* X M Y^*  + \mathcal P_X Z^* X M X^* Y X^{\dagger}+ \mathcal P_X Z^* X M X^* Z \mathcal P_X 
+\mathcal P_XK\mathcal P_X \nonumber \\
&= BB^* +\mathcal P_XK\mathcal P_X,
\end{eqnarray}
where $B=YM_1 + (YX^\dagger)^*(XM_1 ) + \mathcal P_XZ^*X M_1$, and $M=M_1M_1^*$ where $M_1$ is the Cholesky factor of $M$. The first term in~\eqref{finlast} is PSD being a matrix of the form $BB^*$ and the second term is PSD because $K \succeq 0$. Thus~\eqref{finlast}
is PSD being the sum of two PSD matrices. Hence $\Delta  \in \mathbb S(X,Y)$. This shows $``\supseteq"$ in~\eqref{starmatrix} and hence completes the proof. 
\eproof

\begin{remark}{\rm
The assumption that $X^*Y+Y^*X$ is invertible 
is not necessary for the existence of a dissipative mapping which  only requires $X^*Y+Y^*X$ to be positive semidefinite; see Theorem~\ref{thm1:mapping}.  
 This assumption is more of a technical need for the second characterization of dissipative mappings in Theorem~\ref{thm:diffcharac}. However, 
there are situations, such as computing the eigenpair backward error for DH matrices, where one has to find dissipative mappings taking $X \in \C^{n,r}$ to $Y \in \C^{n,r}$ such that $X^*Y+Y^*X$ is invertible. 
For example, consider a DH matrix $J-R$ with $Q=I_n$, $\lambda \in \C$ which is not an eigenvalue of $J-R$, 
and $\hat X \in \C^{n,r}$ with $\text{rank}(\hat X)=r$. Then the structured eigenpair backward error for the DH matrix $J-R$, 
for making $\lambda$ an eigenvalue of multiplicity $r$ and the columns of $\hat X$ as the corresponding eigenvectors, is the smallest norm of the perturbations $\Delta_J-\Delta_R$ such that 
$(J-\Delta_J)+(R-\Delta_R)$ is a DH matrix and 
$((J-\Delta_J)+(R-\Delta_R))\hat X=\lambda \hat X$. By Lemma~\ref{lem:norms_eq_ssh}, the latter condition is equivalent to finding a  dissipative mapping $\Delta$ such that $\Delta \hat X= \hat Y$, where $\hat Y=(\lambda I_n-(J-R)) \hat X$. 
Such a dissipative mapping exists if and only if $\hat X^* \hat Y+\hat Y^*\hat X = \hat X^*(\real{(\lambda)}I_n +R )\hat X  \succeq 0$. Hence if $\real{(\lambda)} >0$, then $\hat X^* \hat Y+\hat Y^*\hat X \succ 0$ since $R \succeq 0$.   
}
\end{remark}

A result analogous to Corollary~\ref{eq:vectorcase} for the vector case (when $m=1$), where the conditions on the vectors $x$ and $y$ are simpler than for  Theorem~\ref{thm:diffcharac}, is stated below. 
\begin{corollary}
Let $x,y \in \mathbb{C}^{n}\setminus{\{0\}}$. 
%
Then $\mathbb{S}(x,y) \neq \emptyset $ if and only if $\real{(x^*y)} \geq 0$.
Moreover, if $\real{(x^*y)}> 0$, then
\begin{eqnarray*} \label{starvec}
\mathbb{S}(x,y) = \{ H + \widetilde{H}(K, G , Z) \mid  K,G, Z\in \mathbb{C}^{n,n},K\succeq 0,~G^*=-G\},
\end{eqnarray*}
where $H$ and $\widetilde H$ are respectively defined as
\begin{equation*}\label{eqH}
H :=\frac{1}{2} \left( \frac{yx^*}{\|x\|^2} - \frac{xy^*}{\|x\|^2}  \right)
+ \frac{\alpha}{2} \left\{ yy^* + (x^*y)\frac{yx^*}{\|x\|^2} 
+ (y^*x)\frac{xy^*}{\|x\|^2}  
+(|x^*y|^2)\frac{ xx^*}{\|x\|^4} \right\}
\end{equation*} 
and
\begin{eqnarray*}\label{eqHtil}
\widetilde{H}(K,G,Z):=\frac{1}{2}\Big( \mathcal P_x K \mathcal P_x + \mathcal P_x G \mathcal P_x +( \alpha(y^*x)+1 ) xx^\dagger Z \mathcal P_x 
+( \alpha(x^*y)-1  ) \mathcal P_x Z^*xx^\dagger  \\ 
+\alpha \mathcal P_x Z^*xy^* + \alpha yx^*Z \mathcal P_x + \alpha \mathcal P_x Z^*xx^*Z \mathcal P_x \Big) , 
\end{eqnarray*}
where $\alpha=\frac{1}{2\real{(x^*y)}}$.
\end{corollary}

\begin{remark}{\rm
Along the lines of Section~\ref{sec:real1}, we can also obtain a second family of real dissipative mappings from Theorem~\ref{thm:diffcharac}, where the restrictions on matrix variables $K,G,Z$ are more simplified than Theorem~\ref{cor:realthm1:mapping}. Indeed, let $X,\, Y \in \C^{n,m}$. Then for $\mathcal X=[X,\,\overline X]$ and $\mathcal Y=[Y,\,\overline Y]$, from Theorem~\ref{thm:diffcharac} whenever there exists a complex dissipative mapping $\Delta$
satisfying $\Delta \mathcal X=\mathcal Y$, then there also exists a real
dissipative mapping. In fact, in view of Lemmas~\ref{realaux} and~\ref{real_lemma2}, the matrix $H$ in Theorem~\ref{thm:diffcharac} is real. Moreover, $H + \widetilde H(K,G,Z)$ with real $K,G,Z \in \R^{n,n}$ such that 
$K \succeq 0$ and $G^T=-G$, gives a family of real dissipative mappings from $X$ to $Y$.
}
\end{remark}

\section{Structured distance to asymptotic instability for  DH systems}\label{seac:DH}
In this section, we study the distance to asymptotic instability for  DH systems
$\dot{x}=(J-R)Qx$ defined in~\eqref{eq:defDH1}.
DH systems are always stable, but they are not necessarily asymptotically stable. 
In~\cite{MehMS16},  the authors have obtained
 various structured distances to asymptotic instability for complex DH systems while perturbing $J$, $R$, or $Q$ (only one matrix at a time).  Similarly in~\cite{MehMS17}, authors derived the real distances to asymptotic instability for real DH systems while perturbing only $R$.
The main tool in the computation of stability radii in~\cite{MehMS16,MehMS17}
 is minimal norm solutions to Hermitian, skew-Hermitian and semidefinite mapping problems. To compute such radii in the case of perturbing both
 $J$ and $R$ simultaneously, we will need minimal norm solutions to complex and real dissipative mappings that map one vector to another vector.
 
%
Therefore, in this section, we exploit the minimal norm solutions to the dissipative mapping problem from previous sections  and derive the structured distance to asymptotic instability for DH systems while perturbing both $J$ and $R$ at a time. By following the terminology in~\cite{MehMS16}, we define the unstructured and the structured distances to asymptotic instability for DH systems as follows:

\begin{definition}\label{def:1}
Let $\mathbb F \in \{\R,\C\}$.
Consider a DH system of the form~\eqref{eq:defDH1} and let $B\in \mathbb F^{n,r}$ and $C \in \mathbb F^{q,n}$ be given
matrices. Then we define
\begin{enumerate}
\item the \textit{unstructured stability radius $r_{\mathbb F}(J,R;B,C)$} of system~\eqref{eq:defDH1} with respect to general perturbations to $J$ and $R$
under the restrictions $(B,C)$ by
\begin{eqnarray}\label{eq:defunst}
r_{\mathbb F}(J,R;B,C)=\inf \Big\{ & \hspace{-1.5cm}\sqrt{{\|\Delta_J\|}_F^2+{\|\Delta_R\|}_F^2}:\Delta_J,\Delta_R\in \mathbb F^{r,q},\nonumber \\
&\Lambda\big((J-R)Q+B(\Delta_J-\Delta_R)CQ\big)\cap i\R \neq \emptyset
\Big\};
\end{eqnarray}
\item the \textit{structured stability radius $r_{\mathbb F}^{S_d}(J,R;B)$} of system~\eqref{eq:defDH1} with respect to structure preserving skew-Hermitian perturbations to $J$ and negative semidefinite perturbations to $R$ from the set
\begin{equation}\label{eq:def_setsd}
S_d^{\mathbb F}(J,R;B)=\left\{(\Delta_J,\Delta_R)\in (\mathbb F^{r,r})^2:~\Delta_J^*=\Delta_J, \Delta_R\preceq 0, R+B\Delta_R B^* \succeq 0 \right\};
\end{equation}
by
\begin{eqnarray}\label{eq:defsst}
r_{\mathbb F}^{S_d}(J,R;B)=\inf \Big\{ & \hspace{-.2cm}\sqrt{{\|\Delta_J\|}_F^2+{\|\Delta_R\|}_F^2}:(\Delta_J,\Delta_R)\in S_d^{\mathbb F}(J,R;B),\nonumber \\
&\Lambda\big((J-R)Q+B(\Delta_J-\Delta_R)B^*Q\big)\cap i\R \neq \emptyset
\Big\}.
\end{eqnarray}
\end{enumerate}
\end{definition}
In the complex case ($\mathbb F=\C$), $r_{\C}(J,R;B,C)$ and $r_{\C}^{S_d}(J,R;B)$ are respectively called the complex unstructured and complex structured stability radii. 
Similarly, in the real case ($\mathbb F=\R$),  $r_{\R}(J,R;B,C)$ and $r_{\R}^{S_d}(J,R;B)$ are respectively called the real unstructured and real structured stability radii.

In order to obtain bounds for the structured stability radius $r_{\mathbb F}^{S_d}(J,R;B)$, we also define the structured eigenvalue backward error as follows 
\begin{eqnarray}\label{eq:def_bacerr}
\eta_{\mathbb F}^{S_d}(J,R;B, \lambda)=\inf \Big\{ & \hspace{-1cm}\sqrt{{\|\Delta_J\|}_F^2+{\|\Delta_R\|}_F^2}:(\Delta_J,\Delta_R)\in S_d^{\mathbb F}(J,R;B),\nonumber \\
&\text{det}\left((J-R)Q+B(\Delta_J-\Delta_R)B^*Q-\lambda I_n\right)=0
\Big\},
\end{eqnarray}
where $\lambda \in \C$, $(J-R)Q$ is a DH matrix, and the perturbation set $S_d^{\mathbb F}(J,R;B)$ is as defined in~\eqref{eq:def_setsd}. Consequently, we have the following result. 
\begin{theorem}\label{thm:radiuserr}
Consider an asymptotically stable DH system of the form~\eqref{eq:defDH1}. Then 
\begin{equation}\label{eq:lboundbacerror}
r_{\mathbb F}^{S_d}(J,R;B)= \inf_{w\in \R} \eta_{\mathbb F}^{S_d}(J,R;B,iw).
\end{equation}
\end{theorem}

\subsection{Complex stability radius}
Consider a complex LTI DH system of the form
\begin{equation}\label{def:complexdh}
\dot{x}(t)=(J-R)Qx(t),
\end{equation}
where $J,R,Q \in \C^{n,n}$ such that $J^*=-J$, $R^*=R\succeq 0$, and 
$Q^*=Q \succ 0$. Here, we study the complex stability radii  $r_{\C}(J,R;B,C)$ and $r_{\C}^{S_d}(J,R;B)$.

Inspired by the proof of~\cite[Proposition 2.1]{HinP86}, we obtain  the following formula for the unstructured stability radius $r_{\C}(J,R;B,C)$. 

\begin{theorem}\label{thm:unstr_rad}
Consider an asymptotically stable DH system of the form~\eqref{def:complexdh}. Let $B \in \C^{n,r}$ and $C \in \C^{q,n}$ be given restriction matrices. Then
$r_{\C}(J,R;B,C)$ is finite if and only if $G(w)=CQ\left(iwI_n-(J-R)Q\right)^{-1}B$ is not zero for some $w \in \R$. In the latter case, we have
\[
r_{\C}(J,R;B,C)=\frac{1}{\sqrt{2}} \inf_{w\in \R}\frac{1}{\|G(w)\|}.
\]
\end{theorem}
\proof
In view of~\eqref{eq:defunst}, note that for any $\Delta_J,\Delta_R\in\C^{r,q}$ 
$\Lambda\big((J-R)Q+B(\Delta_J-\Delta_R)CQ\big)\,\cap \,i\R \neq \emptyset$ if and only if $\text{det}\left((J-R)Q+B(\Delta_J-\Delta_R)CQ-iw I_n\right)=0$
for some $w\in \R$ if and only if $\text{det}\left(I_n-(\Delta_J-\Delta_R)CQ(iw I_n-(J-R)Q)^{-1}B\right)=0$ for some $w\in \R$. Using this in~\eqref{eq:defunst}, we have
\begin{eqnarray}\label{eq:2}
r_{\C}(J,R;B,C) =&\hspace{-2.2cm} \inf \Big\{ \sqrt{{\|\Delta_J\|}_F^2+{\|\Delta_R\|}_F^2}:\Delta_J,\Delta_R\in\C^{r,q},\,w\in \R\nonumber \\
&\text{det}\left(I_n-(\Delta_J-\Delta_R)CQ(iw I_n-(J-R)Q)^{-1}B\right)=0
\Big\}\nonumber\\
&\hspace{-.6cm} = \inf_{w \in \R} \inf \Big\{ \sqrt{{\|\Delta_J\|}_F^2+{\|\Delta_R\|}_F^2}:\Delta_J,\Delta_R\in\C^{r,q},\,x\in \C^{r}\setminus\{0\}\nonumber \\
&(\Delta_J-\Delta_R)CQ(iw I_n-(J-R)Q)^{-1}Bx=x
\Big\}.
\end{eqnarray}
In view of Lemma~\ref{lem:norms_eq}, \eqref{eq:2} becomes
\begin{equation}\label{eq:3}
r_{\C}(J,R;B,C)^2 = \inf_{w \in \R} \inf \Big\{ \frac{{\|\Delta\|}_F^2}{2}:\Delta \in\C^{r,q},\,x\in \C^{r}\setminus\{0\},\Delta CQ(iw I_n\text{$-$}(J\text{$-$}R)Q)^{-1}Bx\text{$=$}x\Big\}.  \\
\end{equation}
Note that $r_{\C}(J,R;B,C) \geq 0$. Therefore $r_{\C}(J,R;B,C)$ is finite if and only if  the inner infimum in~\eqref{eq:3} is finite for some $w \in \R$, which is true if and only if  $G(w):=CQ\left((J-R)Q-iwI_n\right)^{-1}B\neq 0$ for some $w \in \R$. Indeed, if $G(w)\neq 0$ for some $w\in \R$ then for any $x \in \C^{r}\setminus\{0\}$, there exists $\tilde \Delta \in \C^{r,q}$ such that $\tilde \Delta G(w)x=x$. This implies, using~\eqref{eq:3},  that $r_{\C}(J,R;B,C)\leq \frac{{\|\tilde \Delta\|}_F^2}{2}$. Using the fact that ${\|\Delta\|}_F \geq \|\Delta\|$ for any $\Delta$ in~\eqref{eq:3} , we have 
\begin{eqnarray}\label{eq:33}
r_{\C}(J,R;B,C)&\geq & \inf_{w \in \R} \inf \Big\{ \frac{{\|\Delta\|}^2}{2}:\Delta \in\C^{r,q},\,x\in \C^{r}\setminus\{0\},\Delta CQ((J-R)Q-iw I_n)^{-1}Bx=x\Big\}\nonumber \\ \label{eq:inside_1}\\
&=& \frac{1}{2}\inf_{w\in \R} \frac{1}{{\|G(w)\|}^2},\nonumber
\end{eqnarray}
where the last equality follows by slightly modifying the proof
of~\cite[Proposition 2.1]{HinP86}. In fact, we have the equality in~\eqref{eq:inside_1} because for any $w \in \R$ the infimum in the right hand side of~\eqref{eq:inside_1}
is attained for a rank one matrix $\Delta$ for which we have $\|\Delta\|=\frac{1}{{\|G(w)\|}}={\|\Delta\|}_F$~\cite{HinP86}.
\eproof

Next, using Theorem~\ref{thm:radiuserr}, 
we derive bounds for $r_{\C}^{S_d}(J,R;B)$ via bounds for the backward error $\eta_{\C}^{S_d}(J,R;B,\lambda)$.

\begin{theorem}\label{thm:maindhresult}
Consider a DH matrix $(J-R)Q$ and let $w \in \R$ be such that $(J-R)Q-iwI_n$ is invertible. Let $B \in \C^{n,r}$ be of full column rank
and define $\Omega=\text{{\rm null}}\left((I_n-BB^\dagger)(iwI_n-(J-R)Q)\right)$. Suppose that $\Omega \neq \emptyset$. Then
\begin{equation}\label{eq:exact_bacerrtem}
\eta_{\C}^{S_d}(J,R;B, iw)^2\geq \inf_{x\in \Omega} \left\{ \frac{2{\|B^\dagger (iwI_n-(J-R)Q)x\|}^2}{{\|B^*Qx\|}^2}-\frac{|x^*Q(iwI_n-(J-R)Q)x|^2}{{\|B^*Qx\|}^4}
\right\}.
\end{equation}
Moreover, let the infimum in the right hand side of~\eqref{eq:exact_bacerrtem} be attained at $\hat x \in \Omega$ and define
$\hat \Delta_R=-\hat x^*QRQ \hat x \frac{(B^*Q\hat x)(B^*Q \hat x)^*}{{\|B^*Q\hat x\|}^4}$. If $R+B \hat \Delta_R B^* \succeq 0$, then equality holds in~\eqref{eq:exact_bacerrtem}. In this case,
 we have
\begin{equation}
\sqrt{2}\,\sigma_{\min}\left(B^\dagger(iwI_n-(J-R)Q)UW^{-*}\right) \geq \eta_{\C}^{S_d}(J,R;B, iw) \geq \sigma_{\min}\left(B^\dagger(iwI_n-(J-R)Q)UW^{-*}\right),
\end{equation}
where the columns of $U$ form an orthonormal basis for $\Omega$ and $W$ is the
Cholesky factor of $U^*QBB^*QU$.
\end{theorem}
\proof
By definition~\eqref{eq:def_bacerr},
\begin{eqnarray}
\eta_{\C}^{S_d}(J,R;B,iw)
=&\inf \Big\{ \sqrt{{\|\Delta_J\|}_F^2+{\|\Delta_R\|}_F^2}:(\Delta_J,\Delta_R)\in S_d^{\C}(J,R;B),x \in \C^n \setminus\{0\},\nonumber \\
&\left((J-R)Q+B(\Delta_J-\Delta_R)B^*Q-iw I_n\right)x=0
\Big\},\nonumber \\
=&\inf \Big\{ \sqrt{{\|\Delta_J\|}_F^2+{\|\Delta_R\|}_F^2}:(\Delta_J,\Delta_R)\in S_d^{\C}(J,R;B),x \in \C^n \setminus\{0\},\nonumber \\
&B(\Delta_J-\Delta_R)B^*Qx=(iw I_n-(J-R)Q)x
\Big\},\nonumber\\
=&\inf \Big\{ \sqrt{{\|\Delta_J\|}_F^2+{\|\Delta_R\|}_F^2}:(\Delta_J,\Delta_R)\in S_d^{\C}(J,R;B),x \in \Omega\setminus \{0\}, \nonumber \\
&(\Delta_J-\Delta_R)B^*Qx=B^\dagger(iw I_n-(J-R)Q)x
\Big\},\label{eq:aux11}\\
=&\inf \Big\{ {\|\Delta\|}_F:\Delta\in \C^{r,r}, \Delta_R=-(\Delta+\Delta^*) \preceq 0,R+B\Delta_RB^* \succeq0,  \nonumber \\
& x\in \Omega \setminus \{0\},\Delta B^*Qx=B^\dagger(iw I_n-(J-R)Q)x
\Big\},\label{eq:aux12}\\
&\hspace {-1.7cm}\geq\hspace {.2cm}\inf \Big\{ {\|\Delta\|}_F:\Delta\in \C^{r,r}, \Delta_R=-(\Delta+\Delta^*) \preceq 0,x \in \Omega \setminus \{0\},  \nonumber \\
&\Delta B^*Qx=B^\dagger(iw I_n-(J-R)Q)x
\Big\},\label{eq:aux13}
\end{eqnarray}
where we have used Lemma~\ref{lem:auxmap} in~\eqref{eq:aux11} and Lemma~\ref{lem:norms_eq_ssh} in~\eqref{eq:aux12}. If $\Omega \neq \emptyset$, then infimum on the right hand side of~\eqref{eq:aux13}
is finite. In fact, from Corollary~\ref{eq:vectorcase} for any $x\in \Omega$ there exists $\Delta \in \C^{r,r}$ such that $\Delta+\Delta^* \succeq 0$ and
$\Delta B^*Qx=B^\dagger(iw I_n-(J-R)Q)x$ if and only if $\real{\left((B^*Qx)^*B^\dagger(iw I_n-(J-R)Q)x\right)} \geq 0$. Clearly for any $x \in \Omega$, we have
\begin{eqnarray*}
\real{\left((B^*Qx)^*B^\dagger(iw I_n-(J-R)Q)x\right)}  &=&\real{\left(x^*QBB^\dagger(iw I_n-(J-R)Q)x\right)} \\
&=&\real{\left(x^*Q(iw I_n-(J-R)Q)x\right)} \quad (\because x\in \Omega)\\
&=&\real{\left(x^*((iwQ-QJQ)+QRQ)x\right)}\\
&=& \real{(x^*QRQx)}\quad  (\because iwQ-QJQ\,\text{ is\, skew-Hermitian})\\
&=& x^*QRQx \geq 0\quad (\because R \succeq 0,~Q\succ 0).
\end{eqnarray*}
Thus using the minimal norm mapping from Corollary~\ref{eq:vectorcase} in~\eqref{eq:aux13}, we obtain
\begin{equation}\label{eq:exact_bacerrinproof}
\eta_{\C}^{S_d}(J,R;B, iw)^2\geq \inf_{x\in \Omega} \left\{ \frac{2{\|B^\dagger (iwI_n-(J-R)Q)x\|}^2}{{\|B^*Qx\|}^2}-\frac{|x^*Q(iwI_n-(J-R)Q)x|^2}{{\|B^*Qx\|}^4}
\right\}.
\end{equation}
This proves~\eqref{eq:exact_bacerrtem}. Now suppose infimum in the right hand side of~\eqref{eq:exact_bacerrinproof} is attained at $\hat x \in \Omega$, then in view of
Corollary~\ref{eq:vectorcase} consider
\begin{eqnarray*}
&\hat \Delta= \frac{B^\dagger (iwI_n-(J-R)Q)\hat x \hat x^*QB}{{\|B^*Q\hat x\|}^2}
+\frac{(B^\dagger(iwI_n-(J-R)Q)\hat x)^*B^*Q\hat x}{{\|B^*Q\hat x\|}^4}(B^*Q\hat x)(B^*Q\hat x)^*\\
&-\frac{B^*Q\hat x(B^\dagger(iwI_n-(J-R)Q)\hat x)^*}{{\|B^*Q\hat x\|}^2}.
\end{eqnarray*}
Then  $\hat \Delta B^*Q \hat x=B^\dagger(iwI_n-(J-R)Q)\hat x$ and $\hat \Delta+\hat \Delta^* \succeq 0$. Set $\hat \Delta_R=-\frac{\hat \Delta+\hat \Delta^* }{2}$
and $\hat \Delta_J=\frac{\hat \Delta-\hat \Delta^* }{2}$. Then $\hat \Delta_R \preceq 0$. Thus, if $R+B\hat \Delta_R B^* \succeq 0$, then
clearly $(\hat \Delta_J,\hat \Delta_R) \in S_d^{\C}(J,R;B)$, and
\[
{\|\hat \Delta_J\|}_F^2 + {\|\hat \Delta_R\|}_F^2 = {\|\hat \Delta\|}_F^2=  \frac{2{\|B^\dagger (iwI_n-(J-R)Q)\hat x\|}^2}{{\|B^*Q\hat x\|}^2}-\frac{|\hat x^*Q(iwI_n-(J-R)Q)\hat x|^2}{{\|B^*Q\hat x\|}^4}.
\]
This show the equality in~\eqref{eq:exact_bacerrinproof} and hence in~\eqref{eq:aux13}, i.e.,
\begin{equation}\label{eq:exact_bacerrinequality}
\eta_{\C}^{S_d}(J,R;B, iw)^2= \inf_{x\in \Omega} \left\{ \frac{2{\|B^\dagger (iwI_n-(J-R)Q)x\|}^2}{{\|B^*Qx\|}^2}-\frac{|x^*Q(iwI_n-(J-R)Q)x|^2}{{\|B^*Qx\|}^4}
\right\}.
\end{equation}
Now suppose $k=\text{dim}(\Omega)$ and let columns of $U\in \C^{n,k}$ forms an orthonormal basis for $\Omega$. Then $x(\neq 0)\in \Omega$ if and only if
$x=U \alpha$ for some $\alpha\in \C^k \setminus \{0\}$. Using this in~\eqref{eq:exact_bacerrinequality}, we obtain
\begin{equation}\label{eq:exact_bacerrinequality11}
\eta_{\C}^{S_d}(J,R;B, iw)^2= \inf_{\alpha\in \C^k \setminus \{0\}} \left\{ \frac{2{\|B^\dagger (iwI_n-(J-R)Q)U \alpha\|}^2}{{\|B^*QU \alpha\|}^2}-\frac{|\alpha^*U^*Q(iwI_n-(J-R)Q)U \alpha|^2}{{\|B^*Q U \alpha\|}^4}
\right\}.
\end{equation}
Note that $B^*QU$ is a full rank matrix because if suppose $B^*QU \alpha =0$ for some $\alpha\in \C^k \setminus \{0\}$. This implies that
\[
0=B \hat \Delta B^*QU\alpha= (iwI_n-(J-R)Q)U\alpha, 
\]
so that $iw \in \Lambda((J-R)Q)$ which is a contradiction. Thus the matrix $U^*QBB^*QU$ is positive definite. Let $U^*QBB^*QU=WW^*$, where $W$ is the unique Cholesky factor of $U^*QBB^*QU$, then by using 
$y=W^* \alpha$ in \eqref{eq:exact_bacerrinequality11}, we have
\begin{eqnarray}\label{eq:exact_bacerrinequality12}
&\hspace{-9.5cm}\eta_{\C}^{S_d}(J,R;B, iw)^2=\nonumber \\
& \inf_{y\in \C^k \setminus \{0\}}
\left\{ \frac{2{\|B^\dagger (iwI_n-(J-R)Q)U W^{-*}y\|}^2}{{\|y\|}^2}-\frac{|y^*W^{-1}U^*Q(iwI_n-(J-R)Q)U W^{-*}y|^2}{{\|y\|}^4}
\right\}.
\end{eqnarray}
In view of the Cauchy-Schwarz inequality, note that for every $y\in \C^k\setminus\{0\}$ we have
\begin{eqnarray*}
\frac{2{\|B^\dagger (iwI_n-(J-R)Q)U W^{-*}y\|}^2}{{\|y\|}^2} &\geq& \frac{2{\|B^\dagger (iwI_n-(J-R)Q)U W^{-*}y\|}^2}{{\|y\|}^2}\\
&&-\frac{|y^*W^{-1}U^*Q(iwI_n-(J-R)Q)U W^{-*}y|^2}{{\|y\|}^4} \\
&\geq&
\frac{{\|B^\dagger (iwI_n-(J-R)Q)U W^{-*}y\|}^2}{{\|y\|}^2}.
\end{eqnarray*}
Using this in~\eqref{eq:exact_bacerrinequality12} yields
\begin{eqnarray*}
2 \inf_{y \in \C^k \setminus \{0\}}\frac{{\|B^\dagger (iwI_n-(J-R)Q)U W^{-*}y\|}^2}{{\|y\|}^2} &\geq& \eta_{\C}^{S_d}(J,R;B, iw)^2\\
 &\geq&
 \inf_{y \in \C^k \setminus \{0\}} \frac{{\|B^\dagger (iwI_n-(J-R)Q)U W^{-*}y\|}^2}{{\|y\|}^2},
\end{eqnarray*}
and hence
\begin{equation*}
\sqrt{2}\,\sigma_{\min}\left(B^\dagger(iwI_n-(J-R)Q)UW^{-*}\right) \geq \eta_{\C}^{S_d}(J,R;B, iw) \geq \sigma_{\min}\left(B^\dagger(iwI_n-(J-R)Q)UW^{-*}\right).
\end{equation*}
This completes the proof.
\eproof

By using Theorem~\ref{thm:maindhresult} in Theorem~\ref{thm:radiuserr}, we obtain a lower bound for the structured stability radius $r_{\C}^{S_d}(J,R;B)$ as follows. 

\begin{theorem}\label{thm:explicitcompresult}
Consider an asymptotically stable DH system of the form~\eqref{def:complexdh}. Let $B \in \C^{n,r}$ be of full column rank. For $w \in \R$, define $\Omega_w:=\text{{\rm null}}\left((I_n-BB^\dagger)(iwI_n-(J-R)Q)\right)$. Suppose that $\Omega_w \neq \emptyset$ for some $w\in \R$. Then
\begin{equation*}
r_{\C}^{S_d}(J,R;B)^2\geq \inf_{w\in \R} \inf_{x\in \Omega_w} \left\{ \frac{2{\|B^\dagger (iwI_n-(J-R)Q)x\|}^2}{{\|B^*Qx\|}^2}-\frac{|x^*Q(iwI_n-(J-R)Q)x|^2}{{\|B^*Qx\|}^4}
\right\}.
\end{equation*}
\end{theorem}

\begin{remark}{\rm We conclude the section with a few remarks about Theorem~\ref{thm:maindhresult}.
\begin{itemize}
\item Since $\hat \Delta_R$ in Theorem~\ref{thm:maindhresult} is a rank-one matrix with only one negative eigenvalue, if $R \succ 0$ and $R+B\hat \Delta_R B^*$ is singular, then from~\cite[Lemma 4.4]{MehMS16} $R+B\hat \Delta_R B^*$ is positive semidefinite, and this implies the equality in~\eqref{eq:exact_bacerrtem}.

\item As a particular case of Theorem~\ref{thm:maindhresult}, i.e.,  when $w=0$, we
obtain the distance to singularity $d^{S_d}(J,R;B):=\eta^{S_d}(J,R;B,0)$ with respect to structure-preserving perturbations to both $J$ and $R$ from the set 
$S_d(J,R;B)$ which is an analogous result to~\cite[Theorem 6.2]{MehMS16}. Note that in~\cite{MehMS16}, the structured distances to singularity for DH matrices were obtained for structure-preserving perturbations to $Q$ only.


\end{itemize}
}
\end{remark}

\subsection{Real stability radius}
Consider a real LTI DH system of the form
\begin{equation}\label{def:realdh}
\dot{x}(t)=(J-R)Qx(t),
\end{equation}
where $J,R,Q \in \R^{n,n}$ such that $J^T=-J$, $R^T=R\succeq 0$, and 
$Q^T=Q \succ 0$. In this section, we discuss the real unstructured distance  $r_{\R}(J,R;B,C)$ defined by~\eqref{eq:defunst}  and the real structured distance $r_{\R}^{S_d}(J,R;B)$ defined by~\eqref{eq:defsst} for real DH systems. For this, define, for a given $M \in \C^{p,m}$, 
\begin{equation} \label{eq:muRpM} 
\mu_{\R,q}(M):=\left(\inf \left\{{\|\Delta\|}_{q}:~\Delta \in \R^{m,p},~\text{det}(I_m-\Delta M)=0\right\}\right)^{-1},
\end{equation}
where $q \in \{2,F\}$. 
Let us first state a result from~\cite{QiuBRDYD95} that will be useful in determining the unstructured distance 
$r_{\R}(J,R;B,C)$. 
\begin{theorem}[\cite{QiuBRDYD95}]\label{thm:real11}
We have 
\begin{equation*}
\mu_{\R,2}(M)=\inf_{\gamma \in (0,1]} \sigma_2\left(\mat{cc} \real{M} & -\gamma \imag{M}\\ \gamma^{-1} \imag{M} & \real{M}
\rix
\right),
\end{equation*}
where $\sigma_2(A)$ is the second largest singular value of a matrix $A$. Furthermore, an optimal $\Delta$ that attains the value of $\mu_{\R,2}(M)$
can be chosen of rank at most two. 
\end{theorem}

Let us now prove a result bounding $\mu_{\R,F}(M)$ using $\mu_{\R,2}(M)$. 
\begin{theorem} \label{thm:real12}
We have  
\begin{equation}\label{eq:1realthm}
\frac{1}{\sqrt{2}}\mu_{\R,2}(M) \leq \mu_{\R,F}(M) \leq \mu_{\R,2}(M).
\end{equation} 
\end{theorem} 
\proof 
By definition,
\begin{eqnarray}\label{eq:2realthm}
\mu_{\R,2}(M)^{-1}&=&\inf \left\{{\|\Delta\|}_{2}:~\Delta \in \R^{m,p},~\text{det}(I_m-\Delta M)=0\right\}\nonumber\\
&\leq &\inf \left\{{\|\Delta\|}_{F}:~\Delta \in \R^{m,p},~\text{det}(I_m-\Delta M)=0\right\}\nonumber\\
&=& \mu_{\R,F}(M)^{-1},
\end{eqnarray}
since for any $\Delta \in \R^{m,p}$, we have ${\|\Delta\|}_{2} \leq {\|\Delta\|}_{F}$. Now suppose $\hat \Delta  \in \R^{m,p} $ such that $\mu_{\R,2}(M)^{-1}={\|\hat \Delta\|}_{2} $. Then from the first part $\hat \Delta$ is of rank at most two. This implies that  $\text{det}(I_m-\hat \Delta M)=0$ and ${\|\hat \Delta\|}_{F} \leq \sqrt{2} {\|\hat \Delta\|}_{2}$, which yields that 
\begin{equation}\label{eq:3realthm}
\mu_{\R,F}(M)^{-1} \leq  \sqrt{2} {\|\hat \Delta\|}_{2} = \mu_{\R,2}(M).
\end{equation}
Hence from~\eqref{eq:2realthm} and~\eqref{eq:3realthm}, 
we obtain~\eqref{eq:1realthm}.
\eproof

As an application of the above theorem we bound the real unstructured distance $r_{\R}(J,R;B,C)$.

\begin{theorem}\label{thm:realunstrmain}
Consider a real asymptotically stable DH system of the form~\eqref{def:realdh} and let $B \in \R^{n,r}$ and $C \in \R^{q,n}$ be restriction matrices. Then 
\[
\frac{1}{\sqrt{2}}\,\inf_{w\in \R} \mu_{\R,2}(G(w))^{-1} \leq  r_{\R}(J,R;B,C) \leq  \inf_{w \in \R} \mu_{\R,2}(G(w))^{-1},
\]
where
$\mu_{\R,2}(G(w))^{-1}= \inf_{\gamma \in (0,1]} \sigma_2\left(\mat{cc} \real{G(w)} & -\gamma \imag{G(w)}\\ \gamma^{-1} \imag{G(w)} & \real{G(w)} \rix \right)$ and  $G(w)=CQ\left(iwI_n-(J-R)Q\right)^{-1}B$.
\end{theorem}
\proof
By Definition~\eqref{eq:defunst}, we have 
\begin{eqnarray}
& \hspace{-11cm} r_{\R}(J,R;B,C) \nonumber \\
&= \inf \Big\{\sqrt{{\|\Delta_J\|}_F^2+{\|\Delta_R\|}_F^2}:\Delta_J,\Delta_R\in\R^{r,q},\Lambda\big((J-R)Q+B(\Delta_J-\Delta_R)CQ\big)\cap i\R \neq \emptyset \Big\} \nonumber \\
&\hspace{-5.8cm} = \inf \Big\{ \sqrt{{\|\Delta_J\|}_F^2+{\|\Delta_R\|}_F^2}:\Delta_J,\Delta_R\in\R^{r,q},\,w\in \R\nonumber \\
&\text{det}\left((J-R)Q+B(\Delta_J-\Delta_R)CQ-iw I_n\right)=0
\Big\} \nonumber \\
&\hspace{-5.7cm} = \inf_{w\in \R}\inf \Big\{ \sqrt{{\|\Delta_J\|}_F^2+{\|\Delta_R\|}_F^2}:\Delta_J,\Delta_R\in\R^{r,q},\nonumber \\
&\text{det}\left(I_n-(\Delta_J-\Delta_R)CQ(iw I_n-(J-R)Q)^{-1}B\right)=0
\Big\} \nonumber \\
&\hspace{-1.8cm} =\inf_{w\in \R}\inf \Big\{ \frac{{\|\Delta\|}_F}{\sqrt{2}}:\Delta\in\R^{r,q},
\text{det}\left(I_n-\Delta CQ(iw I_n-(J-R)Q)^{-1}B\right)=0
\Big\} \label{eq1:proofunstr} \\
&\hspace{-9.2cm} =\frac{1}{\sqrt{2}} \inf_{w\in \R} \left(\mu_{\R,F}G(w)\right )^{-1},\label{eq2:proofunstr}
\end{eqnarray}
where $G(w)=CQ\left(iwI_n-(J-R)Q\right)^{-1}B$. 
The equality~\eqref{eq1:proofunstr} follows from Lemma~\ref{lem:norms_eq}, and~\eqref{eq2:proofunstr} follows by definition of $\mu_{\R,F}G(w)$ in~\eqref{eq:muRpM}. Thus the result follows immediately using the inequality~\eqref{eq:3realthm} in Theorem~\ref{thm:real12}.
\eproof

We note that a result similar to Theorem~\ref{thm:maindhresult} for the real backward error, $\eta_{\R}^{S_d}(J,R;B)$, defined in~\eqref{eq:def_bacerr}  
may be obtained by using minimal-norm real dissipative mappings from Theorem~\ref{cor:realthm1:mapping}. Thus in view of Theorem~\ref{thm:radiuserr} we can obtain a lower bound for the real structured stability radius $r_{\R}^{S_d}(J,R;B)$.
In the following we state this result for $\eta_{\R}^{S_d}(J,R;B)$ and skip its proof as it is similar to the proof of Theorem~\ref{thm:maindhresult}.

\begin{theorem}\label{thm:realstrdist}
Consider a real DH system of the form~\eqref{def:realdh}. Let $w \in \R$ be such that \mbox{$iwI_n-(J-R)Q$} is invertible. Let $B \in \R^{n,r}$ be of full column rank and let $\Omega$ be the subset of $\text{null}\left((I_n-BB^\dagger)(iwI_n-(J-R)Q) \right)$ such that if $x \in \Omega$ then
$x^*QRQx \geq |x^T(iwQ+QRQ)x|$. Suppose that $\Omega \neq \emptyset$. Then 
\begin{equation}\label{eq1:thrmrealstr}
\eta_{\R}^{S_d}(J,R;B,iw)^2 \geq \inf_{x\in \Omega} \left (
2{\left\| \mathcal Y \mathcal X^\dagger\right\|}_F^2-\text{trace}\left((\mathcal Y \mathcal X^\dagger)^*(\mathcal X \mathcal X^\dagger)(\mathcal Y \mathcal X^\dagger)\right)\right),
\end{equation}
where for $x \in \Omega$, $\mathcal Y=\mat{cc} B^\dagger(iwI_n-(J-R)Q)x & \overline{B^\dagger(iwI_n-(J-R)Q)x}\rix$ and $\mathcal X=\mat{cc}B^TQx & \overline{B^TQx} \rix$. Moreover, let the infimum in the right hand side of~\eqref{eq1:thrmrealstr} be attained at $\tilde x \in \Omega$ 
 and define $\tilde \Delta_R :=  \mathcal {\tilde Y} \mathcal {\tilde X}^\dagger
 -{(\mathcal {\tilde Y} \mathcal {\tilde X}^\dagger)}^* \mathcal P_\mathcal {\tilde X}$, where 
$\mathcal {\tilde Y}=\mat{cc} B^\dagger(iwI_n-(J-R)Q) \tilde x & \overline{B^\dagger(iwI_n-(J-R)Q) \tilde x}\rix$ and $\mathcal {\tilde X}=\mat{cc}B^TQ \tilde x & \overline{B^TQ \tilde x} \rix$. If $R+B \tilde \Delta_R B^T \succeq 0$, then equality holds in~\eqref{eq1:thrmrealstr}.
\end{theorem}
\proof In view of Theorem~\ref{cor:realthm1:mapping}, the proof is similar to Theorem~\ref{thm:maindhresult}.
\eproof

We close the section by noting that a result analogous to Theorem~\ref{thm:explicitcompresult} can be obtain for real structured stability radius $r_{\R}^{S_d}(J,R;B)$ by using Theorem~\ref{thm:realstrdist} in Theorem~\ref{thm:radiuserr}.

\section{Numerical experiments}

In this section, we illustrate the significance of our distances obtained in 
Theorems~\ref{thm:unstr_rad} and~\ref{thm:explicitcompresult} and compare them 
with those of~\cite{MehMS16}, where various structured stability radii have been obtained for DH systems while perturbing only $J$ or only $R$ one at a time.
In the following, l.b. stands for ``lower bound", $r_{\C}(R;B,C)$ (resp.\ $r_{\C}(J;B,C)$) denotes the unstructured stability radius while perturbing only $R$ (resp.\ $J$) with restriction matrices $B$ and $C$~\cite[Theorem 3.3]{MehMS16}, $r_\C^\mathcal S(J;B)$ denotes the structured stability radius with respect to structure-preserving perturbations to $J$~\cite[Theorem 5.1]{MehMS16}, and $r_{\C}^{\mathcal S_d}(R;B)$ denotes the structured stability radius with respect to structure-preserving negative semidefinite perturbations to $R$~\cite[Theorem 4.2]{MehMS16}. 
 We note that $r_{\C}(R;B,C)=r_{\C}(J;B,C)$~\cite[Theorem 3.3]{MehMS16}. It is also clear that 
$r_{\C}(J,R;B,C) \leq  r_{\C}(R;B,C)$ and the minimum of $r_{\C}^\mathcal S(J;B)$ and $r_{\C}^{\mathcal S_d}(R;B)$ gives an upper bound to $r_{\C}^{\mathcal S_d}(J,R;B)$.

All the experiments are performed in Matlab Version No. 9.1.0 (R2016b). To compute these values (or lower bounds), we proceed as follows.

\begin{itemize}

\item The formulas for the three radii  $r_{\C}(R;B,B^*)$, $r_{\C}^\mathcal S(J;B)$, and $r_{\C}^{\mathcal S_d}(R;B)$ were obtained in~\cite{MehMS16} in terms of some non-convex optimization problems. In~\cite{MehMS16}, the authors used the function 
{\it fminsearch} in Matlab to solve these optimization problems which only gives a local minimum starting from an initial guess. We instead used  the {\it GlobalSearch} solver in Matlab, with its default parameters which attempts to locate the global solution. However, we can only guarantee that the computed solutions are 
lower bounds to the exact distances.


\item $r_{\C}(J,R;B,B^*)$:~a formula for $r_{\C}(J,R;B,B^*)$ is obtained in
Theorem~\ref{thm:unstr_rad} which is $\frac{1}{\sqrt{2}}$ times the unstructured distance 
 $r_{\C}(R;B,B^*)$. 
 

\item $r_{\C}^{\mathcal S_d}(J,R;B)$: a lower bound for 
$r_{\C}^{S_d}(J,R;B)$ is obtained in Theorem~\ref{thm:explicitcompresult}.
We again used the {\it GlobalSearch} solver in Matlab for this to get a good approximation for the lower bound. We note that solving the optimization problem involved in~\eqref{eq:exact_bacerrtem} is challenging, 
and beyond the scope of this 
this paper. 
This paper aims to solve the dissipative mapping problem and show it is useful in engineering applications such as DH systems. 
A  possible future work would be to develop more sophisticated ways to solve~\eqref{eq:exact_bacerrtem}. 

\end{itemize}


\begin{example}{\rm \cite[Example 7.1]{MehMS16}}~{\rm 
Consider a prototype example of disk brake squeal problem with the matrices 
$G,M,K, D$ from~\cite{MehMS16}. We consider the DH system $\dot{x}=(J-R)Qx$ with
\[
J=\mat{cc}G & K+\frac{1}{2}N \\ -(K+\frac{1}{2}N^*) & 0  \rix, \quad
R= \mat{cc}D & 0 \\ 0 & 0 \rix, \quad \text{and}\quad
Q=\mat{cc} M & 0\\0& K\rix^{-1}.
\]
We compute the various distances for the restriction matrices 
$B=I_4$ and $C=I_4$, which are given in the following table: 
\begin{center}
\begin{tabular}{|c|c|c|c|c|}
\hline
 $r_{\C}(J,R;B,B^*)$ & $r_{\C}(R;B,B^*)$ &l.b. to $r_{\C}^{\mathcal S_d}(J,R;B)$ &$r_{\C}^\mathcal S(J;B)$ & $r_{\C}^{\mathcal S_d}(R;B)$  \\ \hline
0.0218   & 0.0308  &  0.0310 &   2.4725 &   5.6149 \\ \hline
\end{tabular}
\end{center}

The lower bound $0.0310 \leq r_{\C}^{\mathcal S_d}(J,R;B)$ is obtained by Theorem~\ref{thm:explicitcompresult}, which as expected shows that the unstructured radius, $r_{\C}(J,R;B,B^*) = 0.0218$, is smaller than the structured radius $r_{\C}^{\mathcal S_d}(J,R;B)$.
If we replace the $J$ matrix by $J_1=\mat{cc}G & K \\ -K & 0  \rix$, then the corresponding results are as follows:
\begin{center}
\begin{tabular}{|c|c|c|c|c|}
\hline
 $r_{\C}(J_1,R;B,B^*)$ & $r_{\C}(R;B,B^*)$ &l.b. to $r_{\C}^{\mathcal S_d}(J_1,R;B)$ &$r_{\C}^\mathcal S(J_1;B)$ & $r_{\C}^{\mathcal S_d}(R;B)$  \\ \hline
0.0826  &  0.1169 &   0.1185   & 2.2340    & 5.7971 \\ \hline
\end{tabular}
\end{center}
In this case, the difference between the unstructured radius $r_{\C}(J_1,R;B,B^*)$ and the structured radius $r_{\C}^{\mathcal S_d}(J,R;B)$ is more significant. This implies that if the magnitude of the structured perturbations to both $J$ and $R$ in Frobenius norm is smaller than $0.1185$ than the system remains stable. 

\begin{example}{\rm

To emphasize more on the stability radii $r_{\C}(J,R;B,B^*)$ and $r_{\C}^{\mathcal S_d}(J,R;B)$ obtained in this paper  and to show that these general distances are indeed different than the ones in~\cite{MehMS16},
%
%
we generate  matrices $J$, $R$, $Q$, $B \in \C^{n,n}$ of different sizes $(n \leq 9)$ randomly following a normal distribution with mean $0$ and standard deviation $1$ (randn$(m,n)$ in Matlab) which we project on the feasible set, that is, $J=-J^*$ and $R,Q \succeq 0$, so that $\dot{x}=(J-R)Qx$ is a DH system, and all restricted stability radii were finite.
The restriction matrices $B$ and $C=B^*$ are chosen to be of full rank. We compute our distances $r_{\C}(J,R;B,B^*)$ and $r_{\C}^{\mathcal S_d}(J,R;B)$ and compare these results with those of~\cite{MehMS16} in Table~\ref{tab:radii}. We see that 
(i) the various unstructured and structured distances are quite a different;
(ii) as expected, the stability radii $r_{\C}(J,R;B,B^*)$ and $r^{\mathcal S_d}(J,R;B)$ are significantly smaller than their counterparts from~\cite{MehMS16} of perturbing only one of $J$ and $R$;
(iii) the results indicate that the lower bound to $r_{\C}^{\mathcal S_d}(J,R;B)$ in some cases  $(n=3,8,9)$ is significantly larger than the unstructured distances (columns 2-3 in Table~\ref{tab:radii}) and reasonably close to the upper bound (minimum of the last two columns in Table~\ref{tab:radii}).
}
\end{example}

\begin{table}[h]
\begin{center}
\caption{Various disances to instability for DH matrices} \label{tab:radii}
\begin{tabular}{|c|c|c|c|c|c|}
\hline
size $n$ & $r_{\C}(J,R;B,B^*)$ & $r_{\C}(R;B)$ & l.b. to $r_{\C}^{\mathcal S_d}(J,R;B)$ &$r_{\C}^\mathcal S(J;B)$ & $r_{\C}^{\mathcal S_d}(R;B)$  
\\\hline
3 &    0.1162  &  0.1644 &   0.7353  & 7.8741 &   1.5038 \\ \hline
4&    0.0047 &   0.0067 &   0.5157  &  4.4828 & 620.6015  \\ \hline
5&    0.0550 &   0.0778 &   0.5108 &   1.2273  &  2.5896 \\ \hline
6&    0.0566 &   0.0801 &  0.2764   & 2.1462    & 1.1335 \\ \hline
7&    0.0165  &  0.0234   & 0.2421   & 0.8284  &  4.7189 \\ \hline
 8&   0.0415  &  0.0587 &   0.9596   & 4.3467   & 1.9909 \\ \hline
 9&    0.0632   & 0.0894    &1.3540    &1.7549    &3.6928 \\ \hline
\end{tabular}
\end{center}
\end{table}
}
\end{example}

\section{Conclusion}

In this paper, we have 
derived necessary and sufficient conditions for the existence of  the dissipative mappings taking $X \in \C^{n,k}$ to $Y \in \C^{n,k}$, 
charactered the set of dissipative mappings, 
and 
found the minimal Frobenius norm solution. 
We have then applied these results to DH systems. In  particular, we have used 
dissipative mappings to derive bounds for the structured distance to asymptotic instability for both complex and real DH systems~\eqref{eq:defDH1} when both $J$ and $R$ are subject to perturbations. 

The bounds computed in this paper involve solving a difficult non-convex optimization problem. Possible future works include  more sophisticated methods to compute~\eqref{eq:exact_bacerrtem}, and 
to study other applications of the dissipative mappings.

\section*{Acknowledgement}  

The authors would like to thank the anonymous reviewers for their careful feedback that helped us improve the paper.

\small

\bibliographystyle{siam}
\bibliography{bibliostable}

\end{document}